\documentclass[titlepage,final,12pt]{article}
\usepackage{graphics}
\usepackage{graphicx}
\usepackage{subfigure}
\usepackage{epsfig}
\usepackage{enumerate}
\usepackage{amsfonts}
\usepackage{theorem}
\usepackage{array}
\usepackage[reqno]{amsmath}
\usepackage{color}

\theoremstyle{change}
\newtheorem{proclaim}{PROCLAIM}[section]
\newtheorem{theorem}[proclaim]{Theorem}

\newtheorem{lemma}[proclaim]{Lemma}
\newtheorem{proposition}[proclaim]{Proposition}
\newtheorem{corollary}[proclaim]{Corollary}

\def\state #1. { \noindent{\bf#1.\enspace}}

\newcommand{\comp}{\,{\raise 1pt \hbox{$\scriptstyle\circ$}}\,}

\newcommand{\dist}{\mathop{\rm dist}}
\newcommand{\epi}{\mathop{\rm epi}}
\newcommand{\hypo}{\mathop{\rm hypo}}

\newcommand{\cl}{\mathop{\rm cl}}

\newcommand{\nt}{\mathop{\rm int}}
\newcommand{\lev}{\mathop{\rm lev}}

\newcommand{\uscfcns}{\mathop{\textrm{usc-fcns}}}

\newcommand{\reals}{{I\kern-.35em R}}
\newcommand{\Reals}{\overline{\reals}}

\newcommand{\natnums}{{{\rm l} \kern -.13em {\rm N} }}
\newcommand{\nats}{{I\kern -.35em N}}
\newcommand{\snats}{{I\kern -.29em N}}
\newcommand{\rats}{{Q\kern -.64em \raise 1pt \hbox{$\scriptstyle |$}\;\,}}
\newcommand{\srats}
    {{Q\kern -.56em \raise 1.2pt \hbox{$\scriptscriptstyle /$}\,}}
\newcommand{\ints}{Z\kern -.46em Z}
\newcommand{\ball}{{I\kern -.35em B}}
\newcommand{\Ex}{{I\kern-.35em E}}
\newcommand{\pluss}{\hskip1pt \raise1pt\vbox{\hrule width6pt \vskip1pt \hrule
                    width6pt} \kern-4pt{\lower1pt\hbox{\vrule height6pt
            \kern1pt\vrule height6pt}}\hskip5pt}
\newcommand{\eop}
    {\hfill{$\vcenter{\hrule height1pt \hbox{\vrule width1pt height5pt
     \kern5pt \vrule width1pt} \hrule height1pt}$} \medskip}
\newcommand{\half}
    {{\raisebox{1pt}{$\frac{1}{2}$}}}

\newcommand{\setd}{{ d \kern -.15em l}}
\newcommand{\hatsetd}{ d \hat{\kern -.15em l }}

\renewcommand{\epsilon}{\varepsilon}
\renewcommand{\phi}{\varphi}

\newcommand{\tto}{\;{\lower 1pt \hbox{$\rightarrow$}}\kern -12pt
           \hbox{\raise 2.5pt \hbox{$\rightarrow$}}\;}
\newcommand{\overto}[1]{\,{\raise 0pt\hbox{$\rightarrow$}}\kern -9pt
     \hbox{\lower 3pt \hbox{$\scriptscriptstyle#1$}}\hskip6pt}
\newcommand{\underto}[1]{\,{\lower 1pt\hbox{$\rightarrow$}}\kern -9pt
     \hbox{\raise 4pt \hbox{$\,\scriptscriptstyle#1$}}\hskip7pt}
\newcommand{\bigoverto}[1]{{\raise 0pt\hbox{$\,\longrightarrow$}}\kern -16pt
     \hbox{\lower 3pt \hbox{$\scriptscriptstyle#1$}}\hskip4pt}
\newcommand{\bigunderto}[1]{\,{\lower 1pt\hbox{$\longrightarrow$}}\kern -16pt
     \hbox{\raise 4pt \hbox{$\,\scriptscriptstyle#1$}}\hskip6pt}
\newcommand{\bigbigto}[2]{\,{\raise 0pt\hbox{$\,\longrightarrow$}}\kern -16pt
     \hbox{\lower 3pt \hbox{$\scriptscriptstyle#2$}}\kern -10pt
     \hbox{\raise 4pt \hbox{$\,\scriptscriptstyle#1$}}\hskip7pt}
\newcommand{\downto}{{\raise 1pt \hbox{$\scriptstyle \,\searrow\,$}}}
\newcommand{\upto}{{\raise 1pt \hbox{$\scriptstyle \,\nearrow\,$}}}

\newcommand{\notimply}
    {\quad\hbox{$\Longrightarrow \kern -14pt {/}$}\hskip6pt\quad}

\newcommand{\Lim}{\mathop{\rm Lim}}

\newcommand{\nLim}{\mathop{\rm Lim}\nolimits}

\newcommand{\low}[1]{{\lower1pt \hbox{$\scriptstyle #1$}}}
\newcommand{\loww}[1]{{\lower2pt \hbox{$\scriptstyle #1$}}}
\newcommand{\high}[1]{{\raise1pt \hbox{$\scriptstyle #1$}}}

\newcommand{\cA}{{\cal A}}
\newcommand{\cB}{{\cal B}}

\newcommand{\cR}{{\cal R}}

\newcommand{\nInnLim}{\mathop{\rm InnLim}\nolimits}
\newcommand{\nOutLim}{\mathop{\rm OutLim}\nolimits}
\newcommand{\nlev}{\mathop{\lev}\nolimits}

\newcommand{\nliminf}{\mathop{\rm liminf}\nolimits}
\newcommand{\nlimsup}{\mathop{\rm limsup}\nolimits}
\newcommand{\ninf}{\mathop{\rm inf}\nolimits}
\newcommand{\nsup}{\mathop{\rm sup}\nolimits}
\newcommand{\nmax}{\mathop{\rm max}\nolimits}

\newcommand{\nargmax}{\mathop{\rm argmax}\nolimits}
\newcommand{\nargmin}{\mathop{\rm argmin}\nolimits}

\newcommand{\bfxi}{\mbox{\boldmath $\xi$}}

\newcommand{\lwdy}[2]{\mathrel{\mathop
        {\raisebox{0.1ex}{\null$#1$}}{\hbox{\kern -1.0em
    {\raisebox{-0.8ex}{$\scriptstyle{\;\to #2}$}}}}}}
\newcommand{\lwwdy}[2]{\mathrel{\mathop
        {\raisebox{0.2ex}{\null$#1$}}{\hbox{\kern -1.0em
    {\raisebox{-1.1ex}{$\scriptstyle{\;\to #2}$}}}}}}
\newcommand{\slwwdy}[2]{\scriptsize{{\mathrel{\mathop
        {\raisebox{0.2ex}{\null$#1$}}{\hbox{\kern -1.0em
    {\raisebox{-1.1ex}{$\scriptstyle{\;\to #2}$}}}}}}}}

\def\eto{\,{\lower 1pt\hbox{$\rightarrow$}}\kern -11pt
     \hbox{\raise 4pt \hbox{$\, \scriptstyle e$}}\hskip7pt}
\def\hto{\,{\lower 1pt\hbox{$\rightarrow$}}\kern -11pt
     \hbox{\raise 4pt \hbox{$\, \scriptstyle h$}}\hskip7pt}
\def\pto{\,{\lower 1pt\hbox{$\rightarrow$}}\kern -11pt
     \hbox{\raise 4.5pt \hbox{$\, \scriptstyle p$}}\hskip7pt}
\def\cto{\,{\lower 1pt\hbox{$\rightarrow$}}\kern -11pt
     \hbox{\raise 4pt \hbox{$\, \scriptstyle c$}}\hskip7pt}

\def\dy#1\until#2{\mathrel{\mathop
        {\null#1}\limits^{\hbox{\lower 1.3ex \hbox{$\scriptstyle{\;\to}$}}}}
        {\hbox{{\kern -.2em${\null}
       ^{\hbox{\raise .2ex \hbox{$\scriptstyle{#2}$}}}$}}}}
\def\hidy#1\until#2{\mathrel{\mathop
        {\null#1}\limits^{\hbox{\lower 1.0ex \hbox{$\scriptstyle{\;\to}$}}}}
        {\hbox{{\kern -.2em${\null}
	^{\hbox{\raise 0.9ex \hbox{$\scriptstyle{#2}$}}}$}}}\!}
\def\lody#1\until#2{\mathrel{\mathop
        {\null#1}\limits_{\hbox{\raise 0.5ex \hbox{$\scriptstyle{\to}$}}}}
        {\hbox{{\kern -.25em${\null}
       _{\hbox{\lower 0.7ex \hbox{$\scriptstyle{#2}$}}}$}}}}
\def\lowdy#1\until#2{\mathrel{\mathop
        {\null#1}\limits_{\hbox{\raise 1.0ex \hbox{$\scriptstyle{\to}$}}}}
        {\hbox{{\kern -.25em${\null}
       _{\hbox{\lower 0.9ex \hbox{$\scriptstyle{#2}$}}}$}}}}

\newcommand{\espl}{\textup{\textrm{e-spl}}}

\newcommand{\bcdot}{\,{\raise .2ex \hbox{$\cdot$}}\,}

\begin{document}


\begin{center}
\begin{large}
{\bf Variational Analysis of Constrained M-Estimators}
\smallskip
\end{large}
\vglue 0.7truecm
\begin{tabular}{lcl}
  \begin{large} {\sl Johannes O. Royset 
                                  } \end{large} & \ \ {\phantom{\&}} \ \ &
   \begin{large} {\sl Roger J-B Wets
   				  } \end{large} \\
  \\
  Operations Research Department &&   Department of Mathematics \\
  Naval Postgraduate School  && University of California, Davis \\
  joroyset@nps.edu && rjbwets@ucdavis.edu
\end{tabular}
\end{center}

\vskip 1.3truecm

\noindent {\bf Abstract}. \quad We propose a unified framework for establishing existence of nonparametric
$M$-estimators, computing the corresponding estimates, and proving their
    strong consistency when the class of functions is exceptionally rich.
    In particular, the framework addresses situations where the class of
    functions is complex involving information and assumptions about shape,
    pointwise bounds, location of modes, height at modes, location of
    level-sets, values of moments, size of subgradients, continuity,
    distance to a ``prior'' function, multivariate total positivity, and
    any combination of the above. The class might be engineered to perform
    well in a specific setting even in the presence of little data. The
    framework views the class of functions as a subset of a particular
    metric space of upper semicontinuous functions under the Attouch-Wets
    distance. In addition to allowing a systematic treatment of numerous
    $M$-estimators, the framework yields consistency of plug-in estimators
    of modes of densities, maximizers of regression functions, level-sets of classifiers, and related
    quantities, and also enables computation by means of approximating
    parametric classes. We establish consistency through a one-sided
law of large numbers, here extended to sieves, that relaxes assumptions of
uniform laws, while ensuring global approximations even under model
misspecification.\\

\vskip 0.5truecm

\halign{&\vtop{\parindent=0pt
   \hangindent2.5em\strut#\strut}\cr

{\bf Date}:\quad \ \today \cr}

\baselineskip=15pt

\section{Introduction}

  It is apparent that the class of functions from which nonparametric
    $M$-estimators are selected should incorporate non-data information about
    the stochastic phenomenon under consideration and also modeling
    assumptions the statistician would like to explore.   In applications, the class can become complex involving shape
    restrictions, bounds on moments, slopes, modes, and supports, limits on
    tail characteristics, constraints on the distance to a ``prior''
    distribution, and so on. The class might be engineered to perform
    well in a particular setting; statistical learning is often carried out
    with highly engineered estimators. An ability to consider rich classes of functions leads to
    novel estimators that even in the presence of relatively little data can produce
    reasonable results.

  Numerous theoretical and practical challenges arise when considering
    $M$-estimators selected from rich classes of functions on
    $\reals^d$, which may even be misspecified, as we need to analyze and
    solve infinite-dimensional random optimization problems with nontrivial
    constraints.  In this article, we leverage and extend results from Variational
    Analysis to build a unified framework for establishing existence of
    such {\it constrained M-estimators}, computing the corresponding
    estimates, and proving their strong consistency. We also show strong
    consistency of plug-in estimators of modes of densities, maximizers of
    regression functions, level-sets of classifiers, and related quantities that likewise account for a variety of constraints.
  In contrast to ``classical'' analysis, Variational Analysis centers
    on functions that abruptly change due to constraints and other sources
    of nonsmoothness and therefore emerges as a natural tool for examining
    $M$-estimators selected from rich classes of functions.

\subsection{Setting and Challenges}

Given $d_0$-dimensional random vectors $X^1$, $X^2$, $\dots$, $X^n$, we
consider constrained $M$-estimators of the form
\begin{equation}\label{eqn:estimatorProblem}
\hat f^n \in \epsilon^n\mbox{-}\nargmin_{f\in F^n} \frac{1}{n} \sum_{j=1}^n \psi(X^j,f) + \pi^n(f),
\end{equation}
where $F^n$ is a class of candidate functions on $\reals^d$, or a subset
thereof, possibly varying with $n$ (sieved), $\psi$ is a loss function such as $\psi(x,f) = -\log f(x)$ (maximum likelihood (ML) estimation of
densities) and $\psi((x,y),f) = (y-f(x))^2$ (least-squares (LS) regression),
$\pi^n$ is a penalty function possibly introduced for the purpose of
smoothing and regularization, and the inclusion of $\epsilon^n\geq 0$
indicates that near-minimizers are permitted. We focus on the iid case, but extensions to non-iid samples is possible within our framework.

The Grenander estimator, the ML estimator over log-concave densities, and the
LS regression function under convexity, just to mention a few constrained
$M$-estimators, certainly exist. However, existence is not automatic. For
rich classes of functions, it is rather common to have an empty set of
minimizers in (\ref{eqn:estimatorProblem}); Section 2 furnishes examples. The
extensive literature on $M$-estimators establishes consistency under rather
general conditions (see, e.g., \cite[Thm. 3.2.2, Cor.
3.2.3]{vanderVaartWellner.96}, \cite[Thm. 5.7]{vanderVaart.98}, and
\cite[Thms. 4.3, 4.8]{vandeGeer.00}). Standard arguments pass through uniform convergence of $n^{-1}\sum_{j=1}^n \psi(X^j,\cdot)$ to
$\Ex[\psi(X^1,\cdot)]$, almost surely or in probability, on a sufficiently large class of functions, which in
turn reduces to checking integrability and total boundedness of the class
under an appropriate (pseudo-)metric, the latter being equivalent to finite
metric entropy. It has long been recognized that uniform convergence is
unnecessarily strong; already Wald \cite{Wald.49} adopted a weaker one-sided
condition. In the central case of ML estimation of densities, an upper bound
on $\psi(x,f) = -\log f(x)$ may not be available and typically force
reformulations in terms of $\psi(x,f) = -\log (f(x) + f^0(x))/2f^0(x)$  and
similar expressions, where $f^0$ is some reference density. Uniform
convergence also gives rise to measurability issues, which may require
statements in terms of outer measures \cite{vanderVaartWellner.96}.

In the presence of rich classes of functions, it becomes nontrivial to
compute estimates as there are no general algorithm for
(\ref{eqn:estimatorProblem}). Approximations in terms of basis functions  are not easily
constructed because the class of functions may neither be a linear space nor
a convex set.

\subsection{Contributions}

In this article, we address the challenges of existence, consistency, and
computations of constrained $M$-estimators by viewing the class of functions
under consideration as a subset of a particular metric space of
upper-semicontinuous (usc) functions equipped with the Attouch-Wets
(aw)-distance\footnote{The aw-distance quantifies distances between sets, in
this case hypographs (also called subgraphs) and the name hypo-distance is
sometimes used; see Sec. 3.}. Although viewing $M$-estimators as minimizers
of empirical processes indexed by a metric space is standard, our {\it
particular} choice is novel. The only precursors are
\cite{RoysetWets.14,RoysetWets.15b}, which hint to developments in this
direction without a systematic treatment. Three main advantages emerge from
the choice of metric space: (i) A unified and disciplined approach to rich
classes of functions becomes possible as the aw-distance can be used across
$M$-estimators. (ii) Consistency of plug-in estimators of modes of densities,
maximizers of regression functions, level-sets of classifiers, and related quantities follows
immediately from consistency of the underlying estimators. (iii) Computation
of estimates becomes viable because usc functions, even when defined on
unbounded sets, can be approximated by certain parametric classes to an
arbitrary level of accuracy in the aw-distance. Moreover, the unified
treatment of rich classes of functions allows for a majority of algorithmic
components to be transferred from one $M$-estimator to another.

We bypass uniform laws of large numbers (LLN) and accompanying metric entropy
calculations, and instead rely on a one-sided {\it lsc-LLN} for which upper
bounds on the loss function $\psi$ becomes superfluous. Thus, concern about
density values near zero and the need for reformulations in ML estimation
vanish. Challenges related to measurability reduces to simple checks on the
loss function that can be stated in elementary terms. Already Wald
\cite{Wald.49} and Huber \cite{Huber.67} recognized the one-sided nature of
(\ref{eqn:estimatorProblem}) and this perspective was subsequently formalized
and refined under the name {\it epi-convergence}; see
\cite{DupacovaWets.88,Wang.96,RoysetWets.15c} for results in the parametric
case and also \cite[Ch. 7]{VaAn}. In the nonparametric case, the use of
epi-convergence to establish consistency of $M$-estimators appears to be
limited to \cite{DongWets.00}, which considers ML estimators of densities
that are selected from closed sets in some separable Hilbert space. Moreover,
either the support of the densities are bounded and the Hilbert space is a
reproducing kernel space or all densities are uniformly bounded from above
and away from zero. Sieves are not permitted. The Hilbert space setting is
problematic as one cannot rely on (strong) compactness to ensure existence of
estimators and their cluster points, and weak compactness essentially limits
the scope to convex classes of functions. In addition to going much beyond ML
estimation, our particular choice of metric space addresses issues about
existence. We also provide a novel consistency result that extends the reach
of the lsc-LLN to sieves, which is of independent interest in optimization
theory.

Without insisting on uniform approximations, the lsc-LLN
establishes convergence in some sense across the whole class of functions.
Thus, consistency results are not hampered by model misspecification or other
circumstances under which an estimator is constrained away from an actual
(true) function. They only need to be interpreted appropriately, for example
in terms of minimization of Kullback-Leibler divergence. It also becomes
immaterial whether the estimator and the actual function are unique. Under
misspecification in ML estimation, just to mention one case, there can easily
be an uncountable number of densities that have the same Kullback-Leibler
divergence to the one from which the data is generated. Our results still
hold.

We construct an algorithm for
(\ref{eqn:estimatorProblem}) that under moderate assumptions produces an estimate in a {\it finite} number of iterations if $\epsilon^n>0$
and to converge to an estimate otherwise. The algorithm permits the use of a
wide variety of state-of-the-art optimization subroutines. We demonstrate the
framework in a small study of ML estimation over densities on $[0,1]^2$ that satisfy pointwise upper and lower bounds, have
nonunique modes covering two specific points, are Lipschitz continuous, and
are subject to smoothing penalties.

In our framework, conditions for existence and consistency of estimators
essentially reduce to checking that the class of functions $F^n$ is closed
under the aw-distance. It is well known that the class of concave densities
is closed in this sense. We establish that many other
natural classes of functions are also closed in the aw-distance.
Specifically, we show this for classes defined by convexity, log-concavity,
monotonicity,
     s-concavity, monotone transformations, Lipschitz
     continuity, pointwise upper and lower bounds, location of modes,
     height at modes, location of level-sets, values of moments, size of subgradients, splines, multivariate total positivity of order two,
     and {\it any combination} of the above, possibly under additional assumptions. To the
best of our knowledge, no prior study has established existence and
consistency of $M$-estimators for such a variety of constraints.

We defer the systematic treatment of rates of convergence for $M$-estimators within the proposed framework. Still, because covering numbers of bounded subsets of usc functions under the aw-distance are known \cite{Royset.19}, it is immediately clear that under certain (strong) assumptions rate results can be obtained (see \cite{Royset.19} for preliminary examples), but these are presently not as sharp as those available by means of empirical process theory.

Section 2 provides motivating examples and a small empirical study. Main
results follow in Section 3. Section 4 establishes the closedness of a
variety of function classes under the aw-distance. Section 5 states an
algorithm for (\ref{eqn:estimatorProblem}) and Section 6 gives additional examples. The paper ends with intermediate results and proofs in Section 7.

\section{Motivation and Examples}

The study is motivated, in part, by estimation in the presence of relatively
little data. In such
contexts, constraints in the form of well-selected classes of functions over
which to optimize may become useful. Although statistical
models often aspire to be tuning-free (see for example \cite{CuleSamworthStewart.10a}), models in statistical learning and related application areas are far from being free of tuning \cite{Probst.19}. We
follow that recent trend by considering novel nonparametric estimators
defined by complex constraints, many of which might be tuned to address
specific settings.

\subsection{Role of Constraints}

Analysis using integral-type metrics such as those defined by $L^2$ and
Hellinger distances leads to many of the well-known results for LS regression
and ML estimation of densities. However, difficulties arise with the
introduction of constraints, especially related to closedness and compactness
of the class of function under consideration. For example, consider the class
of bi-constant densities on $[0,1]$, with each density having one value on
$[0,1/2]$ and potentially another value on $(1/2,1]$, that also must satisfy
$f(x)<3/2$ for all $x\in [0,1]$. When the number of samples in $[0,1/2]$ is
sufficiently different from that in $(1/2,1]$, the ML estimator over this
class does not exist as the value of the density in the interval with the
more samples would be pushed up towards the unattainable upper bound. The
break-down is caused by a class of densities that is not closed. Although
rather obvious here, the situation becomes nontrivial in nonparametric cases
involving rich classes of functions that may even be misspecified. In fact,
already the ML estimators over unimodal densities on $\reals$ \cite{Birge.97} and over log-concave densities on $\reals^d$ for $n\leq d$ \cite{DumbgenSamworthSchuhmacher.11} fail to exist.

For another example, suppose that the definition of a class includes the
constraint that the maximizers of the functions should contain a given point
in $\reals^d$. This constraint conveys information or assumption about the
location of modes in the density setting and ``peaks'' in a regression
problem. A sequence of estimates satisfying this constraint may have $L^2$
and Hellinger limits that violate it; the constraint is not closed under
these metrics. Even the simple constraint that $f(\bar x) \geq 1$ for a given
$\bar x\in\reals^d$, which is a constraint on a level-set of $f$, would not
be closed. However,
the constraints on maximizers and such level-sets are indeed closed in the
aw-distance; see Section 2.2 and, more comprehensively, Section 3.4.

Constraints related to maximizers, maxima, and level-sets motivate
the choice of the aw-distance in a profound way as neither pointwise nor
uniform convergence would be satisfactory with regard to those: Pointwise
convergence fails to ensure convergence of maximizers and uniform convergence
applies essentially only to continuous functions defined on compact sets.

\subsection{Example Formulation and Result} As a concrete example of a rich class of densities in ML
estimation on $\reals^d$, suppose that $\alpha,\kappa\geq 0$; $C,D\subset
\reals^d$; $I\subset [0,\infty]$ is closed; $g,h:\reals^d\to [0, \infty)$,
with $h$ being usc and also satisfying $\int h(x) dx <\infty$; and
\begin{align}\label{eqn:FinExample}
  F = &\Big\{f:\reals^d\to [0,\infty]~\Big|~ f \mbox{ usc}, ~\int f(x) dx = 1,\\
      & ~~~~~~~C\subset \nargmax_{x\in\reals^d} f(x),~D\subset \nlev_{\geq\alpha} f, ~\nsup_{x\in \reals^d} f(x) \in I,\nonumber\\
      & ~~~~~~~g(x) \leq f(x)  \leq h(x), ~|f(x) - f(y)|\leq \kappa\|x-y\|_2, \forall x,y\in \reals^d\Big\},\nonumber
\end{align}
where $\nlev_{\geq\alpha} f = \{x\in\reals^d~|~f(x)\geq \alpha\}$ is an upper
level-set of $f$. The second line restricts the consideration to
densities with (global) modes covering $C$ and ``high-probability regions''
covering $D$. Neither $C$ nor $D$ need to be singletons. Although there are
some efforts towards accounting for information about the location of modes
(see for example \cite{DossWellner.16}), the generality of these constraints
is unprecedented. The third line permits nearly arbitrary pointwise bounds. In settings with little data but substantial experience
about what an estimate ``should'' look like, such constraints can be helpful
modeling tools. The last constraint restricts the class to Lipschitz
continuous functions with modulus $\kappa$.

Properties of the ML estimator on this class is stated next. Section 7
furnishes the proof and those of most subsequent results. Let $\nats = \{1, 2,
\dots\}$.

\begin{proposition}\label{prop:example} Suppose that $X^1, X^2, \dots$ are iid random vectors,
each distributed according to a density $f^0:\reals^d \to [0,\infty]$, $F$ in \eqref{eqn:FinExample} is
nonempty, and $\{\epsilon^n\geq 0, n\in\nats\}\to 0$. Then the following hold
almost surely:
\begin{enumerate}[(i)]

\item For all $n\in\nats$, there exists $\hat f^n \in
    \epsilon^n\mbox{-}\nargmin_{f\in F} \{-n^{-1}\sum_{j=1}^n \log f(X^j)\}.$

\item  Every cluster point (under the aw-distance) of $\{\hat f^n, n\in
    \nats\}$, of which there is at least one, minimizes the
    Kullback-Leibler divergence to $f^0$ over the class $F$.

\item If $f^0 \in F$, then $f^0$ is the limit (under the aw-distance)
    of $\{\hat f^n, n\in \nats\}$.

\end{enumerate}
\end{proposition}
The proposition establishes that despite the rather rich class, ML estimators
exists and they are consistent, even under model
misspecification, which is especially relevant in the presence of many constraints. The specific case examined here is only an illustration. Section 3 provides general results along these lines. The framework aims to simplify the analysis of {\it new} estimators constructed by adding and/or removing constraints. As we see below, the analysis reduces largely to checking closedness and nonemptiness.

\subsection{Empirical Results}

\begin{figure}[h!]
\begin{minipage}[t]{0.46\textwidth}
\psfig{file=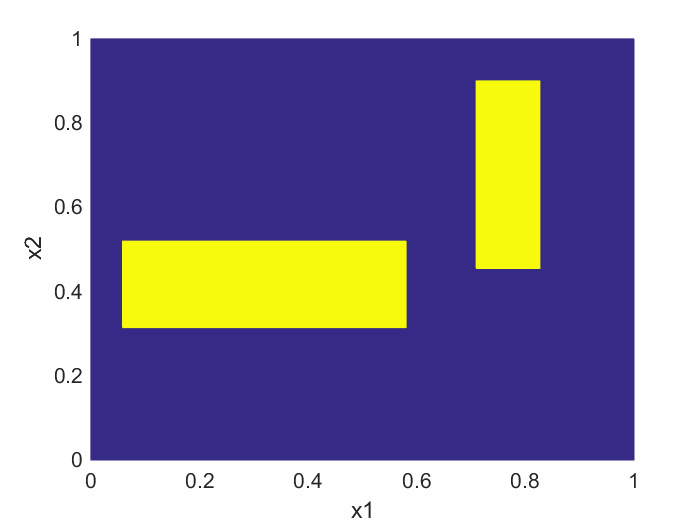,width=1.1\textwidth}
\end{minipage}
\hspace*{\fill} 
\begin{minipage}[t]{0.46\textwidth}
\psfig{file=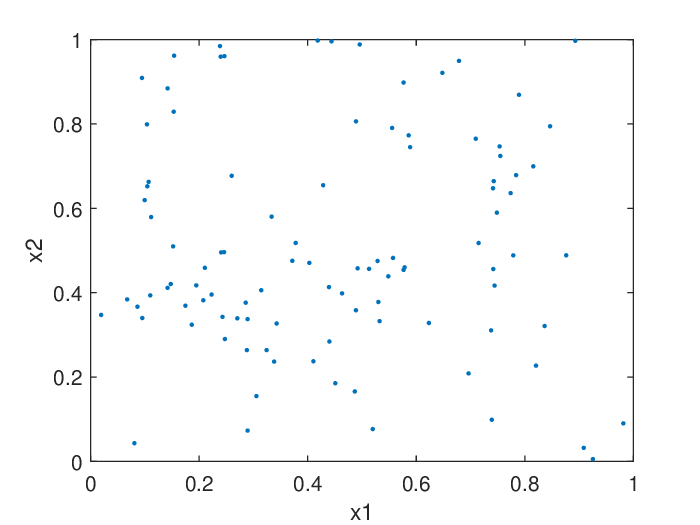,width=1.1\textwidth}
\end{minipage}
\caption[]{Top view of actual density (left) and sample of size $n=100$ (right). }
\label{fig:densities}
\end{figure}

We consider ML estimation of the mixture of three uniform densities on
$[0,1]^2$ depicted in Figure \ref{fig:densities}(left). The resulting mixture
density $f^0$ has height $f^0(x) = 3$ for $x$ in the areas colored yellow and
$f^0(x) = 0.6150$ elsewhere. Using a sample of size 100 shown in Figure
\ref{fig:densities}(right), we compute a penalized ML estimate over the class
of functions
\begin{align*}
F = & \Big\{f:[0,1]^2\to [\alpha,\beta]~ \Big| \int f(x) dx = 1,~\{\bar x, \bar y\} \subset \nargmax_{x\in [0,1]^2} f(x),\\
   & ~~ |f(x) - f(y)|\leq \kappa\|x-y\|_2, \forall x,y\in [0,1]^2,\\
                                  & ~~\mbox{piecewise affine on simplicial complex partition}\Big\}.
\end{align*}
A simplicial complex partition divides $[0,1]^2$ into $N$ equally sized
triangles; see Section 5 for details and the fact that optimization over $F$
can be reduced to solving a finite-dimensional convex problem. As discussed
there, $F$ can be viewed as an approximation, introduced for computational
reasons, of the class obtained from $F$ by relaxing the piecewise affine
restriction. We also adopt the penalty term $\pi(f) = \lambda \sum_{k=1}^N
\|g_i\|_1$, where $g_i$ is the gradient of the $i$th affine function defining
$f$. In the results reported here, $\kappa = 100$ with $\bar x = (0.4702,
0.4657)$ and $\bar y = (0.7746, 0.7773)$. We observe that $F$ is misspecified
as $f^0$ is not Lipschitz continuous.

\begin{figure}[h!]
\begin{minipage}[t]{0.46\textwidth}
\psfig{file=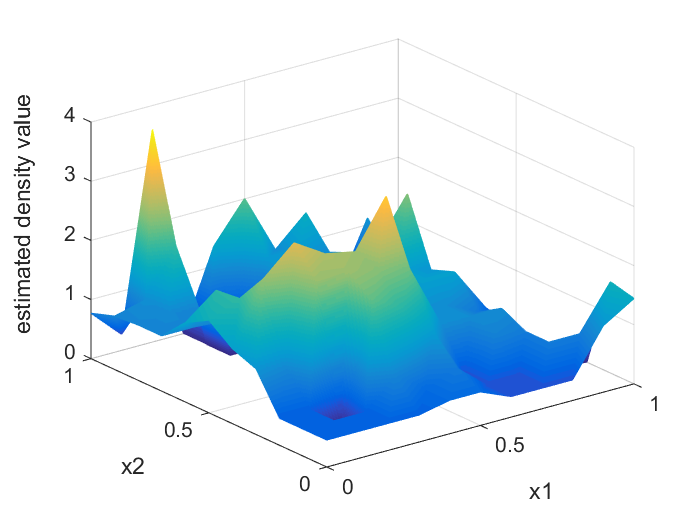,width=1.1\textwidth}
\end{minipage}
\hspace*{\fill} 
\begin{minipage}[t]{0.46\textwidth}
\psfig{file=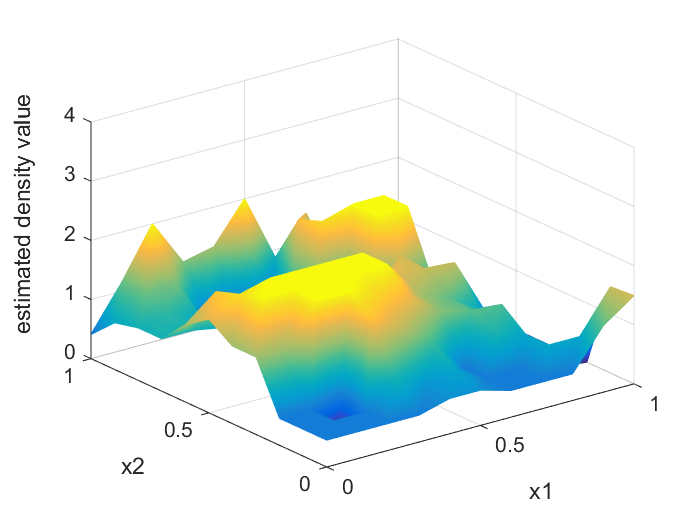,width=1.1\textwidth}
\end{minipage}
\caption[]{Estimates using $n=100$ without (left) and with (right) argmax-constraint. }
\label{fig:densities1}
\end{figure}

Figure \ref{fig:densities1} illustrates the effect of including the
argmax-constraint for the case with $\lambda = 0.05$, $\alpha = 0.0001$,
$\beta = 10000$, and $N = 200$. In the left portion of the figure, the
argmax-constraint is not used and, visually, the errors are large. In the
right portion, the argmax-constraint is included and indications of the
actual density emerges. This and other experiments show that
argmax-constraints regularize the estimates in some sense.

\begin{figure}[h!]
\begin{minipage}[t]{0.46\textwidth}
\psfig{file=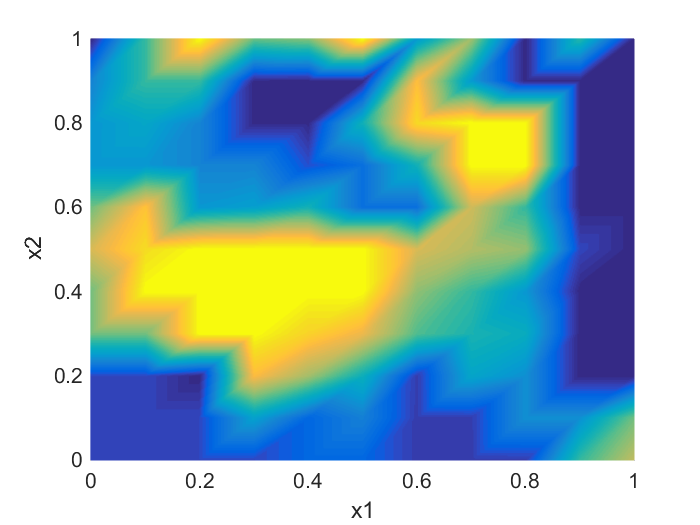,width=1.1\textwidth}
\end{minipage}
\hspace*{\fill} 
\begin{minipage}[t]{0.46\textwidth}
\psfig{file=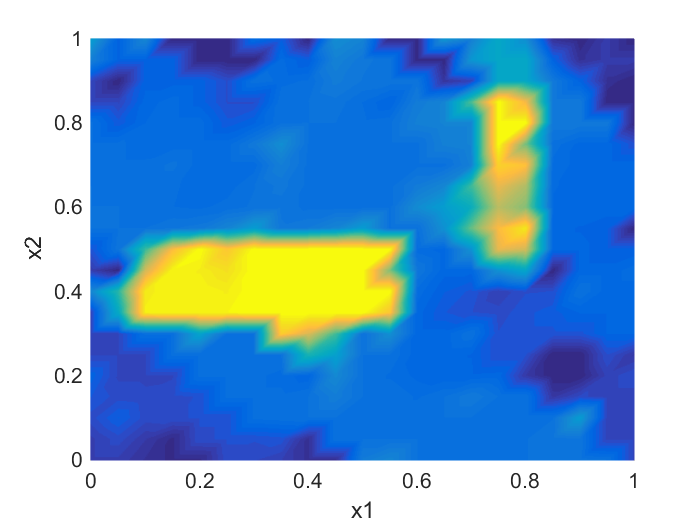,width=1.1\textwidth}
\end{minipage}
\caption[]{Estimates using sample size $n=100$ (left) and $n=1000$ (right). }
\label{fig:densities2}
\end{figure}

If $\alpha$ is increased to 0.3075 and $\beta$ lowered to $4.5$, i.e., 50\%
below and above the lowest and highest point of $f^0$, the estimate with
argmax-constraint is slightly improved; see Figure \ref{fig:densities2}(left)
for a top-view of the resulting density. The estimates are quite insensitive to the choice of $\bar x$ and $\bar y$. Over 25 replications with $\bar x$ randomly selected from the box constituting the left portion of $\nargmax_{x\in [0,1]^2} f^0(x)$ and with $\bar y$ randomly selected from the right box, 22 estimates resemble strongly that in Figure \ref{fig:densities2}(left). The remaining three blur together the two peaks of $f^0$. Still, the KL-divergence between $\hat f^n$ and $f^0$ remains close: the mean across the 25 replications is 0.177 and the standard deviation is 0.005. Naturally, a sample size of $n=1000$
improves the estimates significantly; see Figure \ref{fig:densities2}(right),
where now $\lambda = 0.02$ and $N=800$ are used.

Table \ref{tab:times} summarizes typical computing times on a 2.60GHz laptop
using IPOPT \cite{IPOPT} under varying partition size
$N$, sample size $n$, and penalty parameter; $\alpha$ and $\beta$ are as
before. The solver is not tuned for the specific problem instances and times
can certainly be improved. In most cases, the run times are at most a few
seconds. Interestingly, they are nearly constant in the sample size $n$ as
the size of the optimization problem is independent of $n$; see Section 5.2.
Though, run times grow with partition size $N$. We observe that a piecewise
affine density on a partition of $[0,1]^2$ with size $N$ has $3N$ parameters
that needs to be optimized. Thus, the last row in the table implies
overfitting to some extent. The longer run times with penalties ($\lambda>0$)
are caused by additional optimization variables introduced in implementation
of the nonsmooth penalty term. There are well-known techniques for mitigating
this effect, but they are not explored here. Still, the table indicates the
level of computational complexity for constrained $M$-estimator of this kind.
Section 5 includes further discussion.

\begin{table}
\center
\begin{tabular}{r|rrr|rrr}
\hline
partition & \multicolumn{3}{c|}{without penalty ($\lambda = 0$)} & \multicolumn{3}{c}{with penalty ($\lambda > 0$)}\\
size $N$     &  $n=100$ & $n=1000$ & $n=10000$ &  $n=100$ &  $n=1000$ &  $n=10000$\\
\hline
200          &  0.7     & 0.8&    1.0&  1.0 &   1.0&  1.3\\
800          &  1.7     & 1.7&    1.9& 10.4 &  10.3&  9.6\\
3200         &  6.6     &11.7&   14.8& 38.5 &  29.1& 22.0\\
\hline
\end{tabular}
\center\caption{Computing times in seconds.}\label{tab:times}
\end{table}

\section{Existence and Consistency}

After defining the aw-distance and establishing preliminary properties, this section turns to the main results on existence and consistency of estimators.

\subsection{Attouch-Wets Distance}

\begin{figure}
\centering
    \psfig{file=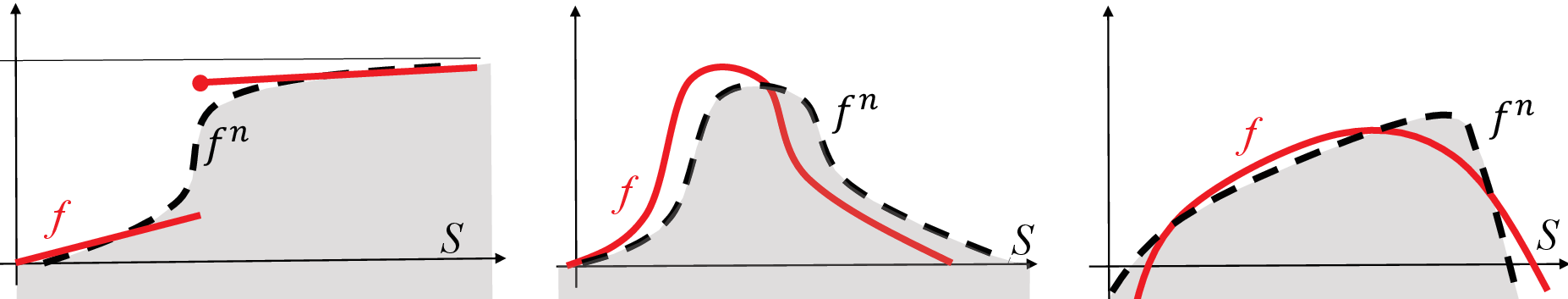,width=0.99\textwidth}
\caption[]{Hypographs of distribution (left), density (middle), and regression functions (right).}
\label{fig:hypoconv}
\end{figure}

Throughout, we consider functions defined on a nonempty and closed set
$S\subset \reals^d$, which may be the whole of $\reals^d$. In the density
setting, $S$ could be thought of as a support. However, we permit densities
to have the value zero, so prior knowledge of the support is not required.
The class $F^n$ in (\ref{eqn:estimatorProblem}) is viewed as a subset of the
(extended real-valued) usc functions on $S$, which is denoted by
 \begin{equation*}
   \uscfcns(S) = \{f:S\to \Reals~|~f \mbox{ usc and } f \not\equiv -\infty\}, \mbox{ with } \Reals = [-\infty,\infty].
 \end{equation*}
Thus, $f\in \uscfcns(S)$ if and only if the {\it hypograph} $\hypo f =
\{(x,\alpha)\in S\times\reals~|~f(x) \geq \alpha\}$ is a nonempty closed
subset of $\reals^d\times \reals$. The class of usc functions is rich enough for
most applications. We equip $\uscfcns(S)$ with the  aw-distance, which
quantifies the distance between hypographs. Figure
\ref{fig:hypoconv} shows $\hypo f^n$ with shading and it appears
``close'' to $\hypo f$. Specifically, let $\dist(z, A)$ be the usual
point-to-set distance between a point $z\in \reals^d\times \reals$ and a
set $A\subset \reals^d\times \reals$; any norm $\|\cdot\|$ can be used. Let
$z^{\rm ctr} \in S\times\reals$. The choice of norm and $z^{\rm
ctr}$ influence the numerical value of the aw-distance, but the resulting
topology on $\uscfcns(S)$ remains unchanged and thus all the stated results
as well. For $f,g\in \uscfcns(S)$, the {\it aw-distance} is defined as
\begin{equation*}
\setd(f,g) = \int_{0}^\infty \setd_\rho(f,g) e^{-\rho} d\rho,
\end{equation*}
where, for $\rho\geq 0$,
\begin{equation*}
\setd_\rho(f,g) = \max \Big\{\big|\dist(z, \hypo f) - \dist(z, \hypo g)\big|~\Big|~\|z-z^{\rm ctr}\| \leq \rho\Big\}.
\end{equation*}
Indeed, $(\uscfcns(S),\setd)$ is a complete separable metric space, for which closed and bounded subsets are compact \cite[Prop. 4.45, Thm. 7.58]{VaAn}. Boundedness can be verified by the inequality $\setd(f,g) \leq 1 + \max\{\dist(z^{\rm ctr},\hypo f)$, $\dist(z^{\rm ctr},\hypo g)\}$ \cite[Prop. 3.1]{Royset.16}. The aw-distance metrizes {\it hypo-convergence}: for $f^n,
f\in\uscfcns(S)$,
\begin{align}\label{eqn:hypoconv}
f^n \mbox{ hypo-converges to } f \Longleftrightarrow & \hypo f^n
\mbox{ set-converges to } \hypo f\nonumber\\
 \Longleftrightarrow & \begin{cases} \forall x^n\to x, ~\nlimsup f^n(x^n) \leq f(x)\\
\forall x~ \exists x^n\to x, ~\nliminf f^n(x^n) \geq f(x)\\
\end{cases}\\
 \Longleftrightarrow & ~\setd(f^n,f)\to 0; \mbox{ simply denoted by } f^n\to f.\nonumber
\end{align}
Set-convergence is in the sense of Painlev\'e{}-Kuratowski\footnote{The {\it outer limit} of a sequence of sets $\{A^n,
n\in\nats\}$ in a topological space, denoted by $\nOutLim A^n$, is the
collection of points to which a subsequence of $\{a^n\in A^n, n\in\nats\}$
converges. The {\it inner limit}, denoted by $\nInnLim A^n$, is the
collection of points to which a sequence $\{a^n\in A^n, n\in\nats\}$
converges. If both limits exist and are equal to $A$, we say that $\{A^n,
n\in\nats\}$ {\it set-converges} to $A$ and write $A^n\to A$ or $\nLim A^n =
A$.}; see \cite[Ch. 7]{VaAn}.

Distribution functions hypo-converge if and only if they converge weakly
\cite{SalinettiWets.86a,SalinettiWets.86,RoysetWets.16b} as illustrated in Figure \ref{fig:hypoconv}(left). Figure \ref{fig:hypoconv}(middle, right) hints
to the fact that modes and maximizers of hypo-converging densities and
regression functions converge to those of limiting functions; see Section
3.4.

In general, $f^n\to f$ does not guarantee pointwise convergence; only $\nlimsup$ $f^n(x)$ $\leq$ $f(x)$ holds for all $x\in S$ by \eqref{eqn:hypoconv}. This issue surfaces in the analysis of (semi)continuity properties of functions on $\uscfcns(S)$. For $\bar x\in S$ and $\rho\geq 0$, let $\ball(\bar x, \rho) = \{x\in S~|~\|\bar x-x\|\leq \rho\}$. We recall that $\{f^n, n\in\nats\}\subset\uscfcns(S)$ is {\it equi-usc} at $\bar x\in S$  when $\nliminf f^n(\bar x)\to \infty$ or when for every $\rho,\epsilon\in (0,\infty)$, there exists $\bar n\in \nats$ and $\delta>0$ such that
\[
\nsup_{x\in \ball(\bar x,\delta)} f^n(x) \leq \max\{f^n(\bar x) + \epsilon, -\rho\} \mbox{ for all } n\geq \bar n.
\]
A class $F\subset\uscfcns(S)$ is equi-usc at $\bar x\in S$ when every sequence $\{f^n\in F, n\in\nats\}$ is equi-usc at $\bar x$. The main consequence of this property is that hypo-convergence implies pointwise convergence \cite[Thm. 7.10]{VaAn}:\\

\begin{proposition}\label{prop:equi-usc_ptwise}{\rm (pointwise convergence).}
If $\{f^n, n\in\nats\}\subset\uscfcns(S)$ is  equi-usc at $\bar x\in S$, then $f^n\to f\in\uscfcns(S)$ implies $f^n(\bar x) \to f(\bar x)$.
\end{proposition}

\vspace{0.3cm}

Although the property is nontrivial, many interesting classes of functions are equi-usc at all, or ``most,'' points in $S$ as seen next. Let $\nt A$ denote the interior of $A\subset\reals^d$.
A log-concave function $f=e^g$ for some concave function $g:S\to\Reals$.

\begin{proposition}\label{prop:equi-usc}{\rm (sufficient conditions for equi-usc).}
Any one of the following conditions suffice for the functions $\{f^n, n\in\nats\}\subset\uscfcns(S)$ to be equi-usc at $\bar x\in S$.

\begin{enumerate}[(i)]

\item The functions are nonnegative, $f^n\to f\in\uscfcns(S)$, and $f(\bar x)=0$.

\item The functions are concave, $f^n\to f\in\uscfcns(S)$, and $\bar x\in \nt \{x\in S~|~f(x)>-\infty\}$.

\item The functions are log-concave, $f^n\to f\in\uscfcns(S)$, and $\bar x \in \nt \{x\in S~|~f(x)>0\}$.

\item The functions are nondecreasing\footnote{Monotonicity of functions on $S$ are always with respect to the partial order induced by inequalities interpreted componentwise, i.e., $f\in\uscfcns(S)$ is nondecreasing (nonincreasing) if $x\leq y$ implies $f^n(x) \leq (\geq) f^n(y)$.} (alternatively, nonincreasing), $f^n\to f\in\uscfcns(S)$, and $f$ is continuous at $\bar x\in \nt S$.

\item The functions are locally Lipschitz continuous at $\bar x$ with common modulus, i.e., there exist $\delta>0$ and $\kappa\in [0,\infty)$ such that $|f^n(x) - f^n(\bar x)|\leq \kappa\|x-\bar x\|$ for all $x\in \ball(\bar x,\delta)$ and $n\in\nats$.

\end{enumerate}
\end{proposition}

Although the aw-distance cannot generally be related to some of the other
common metrics, the
Hellinger and $L^2$
distances tend to zero whenever the aw-distance vanishes under equi-usc and integrability assumptions.

\begin{proposition}\label{prop:integral1}{\rm (connections with other metrics).}
Suppose that $\{f, f^n, n\in\nats\}\subset \uscfcns(S)$, $f^n\to f$, and
for some (measurable) $g:S\to [0,\infty]$, $|f^n(x)|\leq g(x)$ for all $x\in S$
and $n\in\nats$. Then,
\begin{enumerate}[(i)]

\item $L^2_{P}(f^n,f) = \int \big(f^n(x) - f(x)\big)^2 dP(x) \to 0$  provided $\int g^2(x)dP(x)<\infty$ and $\{f^n, n\in\nats\}$ is equi-usc at $P$-a.e. $x\in S$;

\item $H^2(f^n,f) = \half\int(\sqrt{f^n(x)} - \sqrt{f(x)})^2 dx\to 0$
    provided that $f^n\geq 0$, $\int g(x)dx<\infty$, and $\{f^n, n\in\nats\}$ is equi-usc at Lebesgue-a.e. $x\in S$.
\end{enumerate}
\end{proposition}

\subsection{Existence}

Our first main result establishes that existence of an estimator reduces to
having a semi-continuity property for the loss and penalty functions and a
closed and bounded class of functions in the aw-distance.

We recall that a function $\phi:F\to \Reals$ defined on a closed
subset $F$ of $\uscfcns(S)$ is lower-semicontinuous (lsc) if $\nliminf
\phi(f^n) \geq \phi(f)$ for all $f^n\in F\to f$. To clarify earlier
notation\footnote{Throughout we use the common extended real-valued calculus:
$0\cdot \infty = 0$, $\alpha\cdot \infty = \infty$ for $\alpha>0$, $\alpha +
\infty = \infty$ for $\alpha\in \Reals$, and $\alpha -
\infty = -\infty$ for $\alpha\in [-\infty,\infty)$; see \cite[Sec. 1.E]{VaAn}.}, let
$\epsilon\mbox{-}\nargmin_{f\in F} \phi(f)$ $=$ $\{f\in F~|~\phi(f) \leq
\ninf_{g\in F} \phi(g) + \epsilon\}$.

Although our focus is on the existence of $M$-estimators, i.e., minimizers of losses under an empirical distribution, occasionally we consider general distributions and thereby also treat approximation problems. We consider the following general setting \cite[Ch. 14]{VaAn}: For a closed
$F\subset\uscfcns(S)$ and a complete probability space $(S^0, \cB^0, P^0)$, with $S^0\subset \reals^{d_0}$, we
say that $\psi:S^0\times F\to \Reals$ is a {\it random lsc
function} if for all $x\in S^0$, $\psi(x,\cdot)$ is lsc and $\psi$ is
measurable with respect to the product sigma-algebra\footnote{For $F$, we
adopt the Borel sigma-algebra under $\setd$.} on $S^0\times F$.
A random lsc function $\psi:S^0\times F\to \Reals$ is {\it locally
inf-integrable} if for all $f\in F$ there exists $\rho>0$ such
that\footnote{With $\infty - \infty = \infty$, the integral of any measurable function is well-defined. In particular, the present integrand is measurable \cite[Thm. 14.37]{VaAn}, \cite[Prop. 6.3]{DongWets.00}.} $\int \ninf_{g\in F}
\{\psi(x,g)~|~\setd(f,g)\leq \rho\} dP^0(x)
> -\infty$.

\begin{theorem}\label{thm:existence}{\rm (existence of approximation).} Suppose that $\epsilon\geq 0$
and $F$ is a nonempty, closed, and bounded subset of $\uscfcns(S)$ and $(S^0, \cB^0, P^0)$ is a complete probability space.
If $\psi:S^0\times F\to \Reals$ is a locally inf-integrable random lsc
function and $\pi:F\to (-\infty, \infty]$ is lsc, then
\begin{align*}
&\epsilon\mbox{-}\nargmin_{f\in F} \int \psi(x, f) dP^0(x) + \pi(f) \neq\emptyset\\
\mbox{and } & \inf_{f\in F} \int \psi(x, f)dP^0(x) + \pi(f)>-\infty.
\end{align*}
\end{theorem}

\begin{corollary}\label{cor:existence}{\rm (existence of estimator).} Suppose that $\epsilon\geq 0$, $\{x^j\in\reals^{d_0}$, $j=1$, $\dots$, $n\}$,
and $F$ is a nonempty, closed, and bounded subset of $\uscfcns(S)$.
If $\pi:F\to (-\infty, \infty]$ and $\psi(x^j,\cdot):F\to (-\infty, \infty]$
are lsc     for all $j$, then
\[
\epsilon\mbox{-}\nargmin_{f\in F} \frac{1}{n} \sum_{j=1}^n \psi(x^j, f) + \pi(f) \neq\emptyset \mbox{ and } \inf_{f\in F} \frac{1}{n} \sum_{j=1}^n \psi(x^j, f) + \pi(f)>-\infty.
\]
\end{corollary}

For $F$ to be bounded it suffices that there are $x\in S$ and $\alpha\in\reals$ such that for all $f\in F$, $f(x)\geq \alpha$, which becomes trivial for densities and distribution
functions. As we see below, the condition can sometimes be removed.
Many natural classes are closed as indicated in the introduction and
detailed in Section 4. The common penalty function $\pi(f) = \nsup_{x\in S} |f(x)|$ is lsc (cf. Proposition \ref{prop:lscsupnorm}). Familiar loss functions satisfy the lsc requirement
too, at least under certain assumptions. Several examples are furnished including some involving support vector machines (SVM). The class in the next corollary considers concave classifiers in a ``band'' that are also subject to constraints on the location of level-sets.

\begin{corollary}\label{cor:existence-svm}{\rm (existence of concave SVM classifier).} For $g:S\to (-\infty, \infty]$, $h\in \uscfcns(S)$, $\alpha\in\reals$, and $C\subset\reals^d$, suppose that $\{y^j\in \{-1, 1\}, x^j\in \nt S$, $j=1$, $\dots$, $n\}$
and $F=\{f\in\uscfcns(S)~|~ f \mbox{ concave}, ~g(x)\leq f(x) \leq h(x) ~\forall x\in S, ~C\subset \nlev_{\geq \alpha} f\}$. Then, as long as $F$ is nonempty,
\[
\nargmin_{f\in F} \frac{1}{n} \sum_{j=1}^n \max\big\{0, 1 - y^j f(x^j)\big\} \neq\emptyset.
\]
\end{corollary}

When classification errors of different types need to be treated separately, a Neyman-Pearson model leads to the following setting \cite{CannonHowseHushScovel.02,CasasentChen.03,ScottNowak.05}.

\begin{corollary}\label{cor:existence-np}{\rm (existence of robust Neyman-Pearson classifier).} Suppose that $\{x^j\in S$, $j=1$, $\dots$, $n\}$ and $\{z^i\in S$, $i=1$, $\dots$, $m\}$ are associated with $+1$ and $-1$ labels, respectively, and $F$ is a nonempty closed subset of $\uscfcns(S)$.
Then, for open sets $\{Z^i\subset \reals^d, i=1, \dots, m\}$, with $z^i\in Z^i$,
\[
\nargmin_{f\in F} \Big\{\frac{1}{n} \sum_{j=1}^n \max\big\{0, 1 - f(x^j)\big\}~\Big|~ f(z) \leq 0 ~\forall z\in Z^i, i=1, \dots, m \Big\} \neq\emptyset.
\]
\end{corollary}

The corollary establishes the existence of an estimator, defined by a broad class $F$, that minimizes hinge loss across the $+1$ labels and tolerates no training error across the $-1$ labels even after perturbations within sets $Z^i$.

\begin{corollary}\label{cor:exist_ML}{\rm (existence of ML estimator).}
If $F$ is a nonempty closed subset of $\uscfcns(S)$ consisting of nonnegative
functions, $\epsilon\geq 0$, $\{x^j\in S$, $j=1$, $\dots$, $n\}$, and
$f(x^j)<\infty$ for all $j$ and $f\in F$, then
\begin{equation*}
\epsilon\mbox{-}\nargmin_{f \in F} -\frac{1}{n} \sum_{j=1}^n \log f(x^j)\neq\emptyset \mbox{ and } \inf_{f \in F} -\frac{1}{n} \sum_{j=1}^n \log f(x^j)>-\infty.
\end{equation*}
\end{corollary}

We observe that the corollary actually applies to any $f:S\to [0,\infty]$ and
not only densities\footnote{We extend $\alpha\mapsto\log \alpha$ to
$[0,\infty]$ by assigning the end points $-\infty$ and $\infty$,
respectively.}. This fact is beneficial in analysis of estimators for which
the integral-to-one constraint is relaxed, for example, due to computational concerns.
Nevertheless, the constraint enters in many settings and needs a closer examination.
%

If $F$ is the class of
normal densities with mean zero and positive standard deviation, then $F$ is
not closed because there is a sequence in $F$ hypo-converging to a degenerate
density with zero standard deviation. Similarly, if $F$ is the class of
normal densities with standard deviation one, then closedness fails again
since one can construct densities in $F$ hypo-converging to the zero
function. Also classes of bounded densities on a compact set $S$ may not be
closed.  Elimination of such pathological cases is required for a class of densities to be closed. Proposition \ref{prop:example} furnishes a concrete example, while Proposition \ref{prop:integral2} shows that if $F$ is equi-usc at Lebesgue-a.e. $x\in S$ and an integrability condition holds, then $\int f(x) dx = 1$ is closed under hypo-convergence.
The log-concave class exhibits an equi-usc property as established in Proposition \ref{prop:equi-usc}. It is therefore not surprising that the ML estimator over this class exists under a mild condition on the sample \cite{DumbgenSamworthSchuhmacher.11}; see Proposition \ref{prop:exist_ML-logconcave} below.

\begin{corollary}\label{cor:exist_LS}{\rm (LS regression).}
Suppose that $F$ is a nonempty closed subset of $\uscfcns(S)$. If $\{y^j\in \reals, x^j\in S,
j=1, \dots, n\}$ and $F$ is equi-usc at $x^j$, $j=1, \dots, n$, then
\begin{equation*}
\nargmin_{f \in F} \frac{1}{n} \sum_{j=1}^n (y^j - f(x^j))^2 \neq \emptyset.
\end{equation*}
\end{corollary}

Proposition \ref{prop:equi-usc} gives various sufficient conditions for a class of functions to be equi-usc.
The concave functions are equi-usc at ``most'' points according to that proposition and a variant of LS regression that also includes pointwise upper and lower bound, for example introduced to engineer desirable estimates in high-dimensional settings, does indeed exists. This can be established using the same arguments as those supporting Corollary \ref{cor:existence-svm}.

In special cases with relatively simple constraints such as only monotonicity or only convexity, existence of LS estimators are well-known; see \cite{SasabuchiInutsukaKulatunga.83,SeijoSen.11}. The key feature of these special cases is that they reduce in some sense to finite-dimensional problems expressed in terms of the heights $\theta_j = f(x^j), j=1, \dots, n$, and the limit of sequence of such heights generated by feasible functions can easily be shown to be extendable to a feasible function. In the presence of nontrivial constraints that impose restrictions on $f$ at points other than the design points, the situation is more complicated and our systematic approach has merit. In particular, starting from a closed equi-usc class, one can build up closed equi-usc classes through set operations that preserve closedness such as intersections and thereby construct novel estimators that will exist by Corollary \ref{cor:exist_LS}.

\subsection{Consistency}

Our second main result establishes that consistency follows essentially from lower-semicontinuity and one-sided integrability of the loss function and the
closedness of the class under consideration.

\begin{theorem}\label{thm:con}{\rm (consistency).}
Suppose that $X^1, X^2, \dots$ are iid random vectors with values in
$S^0\subset\reals^{d_0}$, $F$ is a closed subset of $\uscfcns(S)$,
$\psi:S^0\times F\to \Reals$ is a locally inf-integrable random lsc
function, and $\pi^n:F\to [0,\infty)$ satisfies $\pi^n(f^n)\to 0$ for every
convergent sequence $\{f^n \in F, n\in \nats\}$. Then, the following hold
almost surely:
\begin{enumerate}[(i)]

\item For all $\{\epsilon^n\geq 0, n\in\nats\}\to 0$,
\begin{align*}
&\nOutLim \Big(\epsilon^n\mbox{-}\nargmin_{f\in F} \frac{1}{n}\sum_{j=1}^n \psi(X^j,f) + \pi^n(f)\Big)\\
&~~~~~~~~~~~~~~~~~~~~~~~~~~~~~~~~~~~~~~~~~~~~~~~~~~ \subset \nargmin_{f\in F} \Ex[\psi(X^1,f)].
\end{align*}

\item There exists $\{\epsilon^n\geq 0, n\in\nats\}\to 0$, such that
\begin{equation*}
\Big(\epsilon^n\mbox{-}\nargmin_{f\in F} \frac{1}{n}\sum_{j=1}^n \psi(X^j,f)+ \pi^n(f)\Big) \to \nargmin_{f\in F} \Ex[\psi(X^1,f)]
\end{equation*}
provided that $\Ex[\psi(X^1,f)]<\infty$ for at least one $f\in F$ and $F$ is bounded.
\end{enumerate}
\end{theorem}

The first conclusion of Theorem \ref{thm:con} guarantees that every cluster
point of sequences constructed from near-minimizers of $n^{-1}\sum_{j=1}^n
\psi(X^j,\cdot) + \pi^n$ is contained in $\nargmin_{f\in F} \Ex[\psi(X^1,f)]$
provided that $\epsilon^n$ vanishes.

Since $\nargmin_{f\in F} \Ex[\psi(X^1,f)]$ may not be a singleton, especially
under model misspecification, there might be a strict inclusion in the first conclusion. For example, let $S=S^0 =
[0,1]$, $F = \{f~|~ f(x) = 1 \mbox{ for } x \in [0,1), f(1)
\in [1,2]\}$, the actual density $f^0$ be uniform on $S$, and $\pi^n(f) =
n^{-1} \nsup_{x\in S} f(x)$. Then, almost surely, $\nargmin_{f\in F} -n^{-1}
\sum_{j=1}^n \log f(X^j)$ $+$ $\pi^n(f)$ $=$ $\{f^0\}$, a strict subset of
$\nargmin_{f\in F} \Ex[-\log f(X^1)]=F$. In this example, the
difficulty is caused by effects on a set of Lebesgue measure zero. However, in more complicated
situations, the concern may be more prevalent. An example is furnished
by the same $f^0$, $S$, and $S^0$, but with $F = \{g^1, g^2\}$, where $g^1(x)
= 1+\delta$ for $x \in [0,1/2]$ and $g^1(x) = 1-\delta$ for $x \in (1/2,1]$,
and $g^2(x) = 1-\delta$ for $x \in [0,1/2]$ and $g^2(x) = 1+\delta$ for $x
\in (1/2,1]$, where $\delta\in (0,1)$, and $\pi^n(f) = n^{-1/2} f(0)$. The
actual density $f^0$ is outside $F$. Then, almost surely, $\nOutLim
\{\nargmin_{f\in F} -n^{-1}\sum_{j=1}^n \log f(X^j)$ $+$ $\pi^n(f)\}$ $=$
$\{g^2\}$, a strict subset of $\nargmin_{f\in F} \Ex[-\log f(X^1)]= F$.

The second conclusion in Theorem \ref{thm:con} guarantees that if
$\epsilon^n$ tends to zero sufficiently slowly, then the inclusion cannot be
strict; near-minimizers of $n^{-1}\sum_{j=1}^n \psi(X^j,\cdot) + \pi^n$
set-converge to $\nargmin_{f\in F} \Ex[\psi(X^1,f)]$. Thus, in this sense,
estimators can converge to {\it any} function in the latter argmin.

A comparison with the common approach to consistency laid out, for example, in \cite[Sec. 3.2.1]{vanderVaartWellner.96} is illuminating. In our notation, \cite[Cor. 3.2.3]{vanderVaartWellner.96} states roughly that if (i) $n^{-1}\sum_{j=1}^n \psi(X^j,f)$ converges in probability to $\Ex[\psi(X^1,f)]$ uniformly in $f$ across $F$, which is permitted to be {\it any} metric space, and (ii) $\Ex[\psi(X^1,\cdot)]$ has a well-separated (unique) minimizer $f^0$ on $F$, then $\hat f^n$ converges in probability to $f^0$. The ability to handle an arbitrary metric space is an advantage over Theorem \ref{thm:con}, but also burdens the user with verifying the well-separability of $f^0$ in the chosen metric. We do not insist on a unique minimizer as discussed above.  The required uniform weak law of large numbers would typically need $\psi(X^1,f)$ to be integrable. In contrast, Theorem \ref{thm:con} insists only on a one-sided integrability condition, which is trivially satisfied when $\psi(x, f)$ is uniformly bounded from below across $x\in S^0$ and $f\in F$ as would be the case for hinge-loss, least-squares, and other common loss functions.

\begin{corollary}{\rm (consistency for concave SVM classifier).}\label{cor:consvm} For $g:\reals^d\to (-\infty, \infty]$, $h\in \uscfcns(\reals^d)$, $\gamma\in\reals$, and $C\subset\reals^d$, suppose that $(X^1,Y^1)$, $(X^2,Y^2), \dots$ are iid random vectors in $\reals^d\times \{-1,1\}$ and $F=\{f\in\uscfcns(\reals^d)$ $|$ $f \mbox{ concave}$, ~$g(x)$ $\leq$ $f(x) \leq h(x) ~\forall x\in \reals^d, ~C\subset \nlev_{\geq \gamma} f\}$.

If $\{\epsilon^n\geq 0, n\in\nats\}\to 0$ and
\begin{equation*} \hat f^n \in
\epsilon^n\mbox{-}\nargmin_{f\in F} \frac{1}{n}\sum_{j=1}^n \max\{0, 1- Y^jf(X^j)\},
\end{equation*}
then, almost surely, $\{\hat f^n, n\in\nats\}$ has at least one cluster point
and every such point $f^\star$ satisfies
\[
f^\star \in \nargmin_{f\in F} \Ex\big[ \max\{0, 1- Y^1f(X^1)\}\big].
\]
Moreover, for a subsequence $\{n_k, k\in\nats\}$ with $\hat f^{n_k}\to f^\star$ and  $\beta<\alpha\in\reals$,
\[
\nOutLim_k \big(\nlev_{\geq \alpha} \hat f^{n_k}\big) \subset \nlev_{\geq \alpha} f^\star \mbox{ and } \nInnLim_k \big(\nlev_{\geq \beta} \hat f^{n_k}\big) \supset \nlev_{\geq \alpha} f^\star.
\]
\end{corollary}

We note that the upper level-sets of $\hat f^n$, which are central in the practical use of the classifier (especially for $\alpha=0$), indeed approximate the ``true'' level-set $\nlev_{\geq \alpha} f^\star$. Without additional assumptions, we are unable to permit $\beta=\alpha$ because it is fundamentally difficult to estimate $\nlev_{\geq \alpha} f^\star$ when $f^\star(x) = \alpha$ on a set of positive measure. For consistency of SVM defined over a subset of a reproducing kernel Hilbert space, we refer to \cite{Steinwart.05}.

The Kullback-Leibler divergence
\[
K(g;f) = \int g(x) \big[\log g(x) - \log f(x)\big] dx
\mbox{ for (measurable) } f,g:S\to [0,\infty]
\]
enters in ML estimation of densities.

\begin{corollary}{\rm (consistency in ML estimation).}\label{cor:conML}
Suppose that $X^1, X^2, \dots$ are iid random vectors, each distributed
according to a density $f^0:S\to [0,\infty]$, $F$ is a closed subset of
$\uscfcns(S)$ with nonnegative functions, and for every $f\in F$ there exists
$\rho>0$ such that $\Ex[\nsup_{g\in F} \{\log g(X^1)~|~\setd(f,g)\leq \rho\}]
<\infty$. If $\{\epsilon^n\geq 0, n\in\nats\}\to 0$ and
\begin{equation*} \hat f^n \in
\epsilon^n\mbox{-}\nargmin_{f\in F} -\frac{1}{n}\sum_{j=1}^n \log  f(X^j),
\end{equation*}
then, almost surely, $\{\hat f^n, n\in\nats\}$ has at least one cluster point
and every such point $f^\star$ satisfies
\[
f^\star \in \nargmin_{f\in F} K(f^0; f).
\]
Under the additional assumption that $F$ contains only densities and $f^0\in
F$, we also have that $f^\star(x) = f^0(x)$ for Lebesgue-a.e. $x\in S$.
\end{corollary}

It is obvious that when there exists an $\alpha\in \reals$ such that $f(x)
\leq \alpha$ for all $f\in F$, then the expectation assumption is satisfied.
In particular, such an $\alpha$ exists if for some $\kappa\in [0,\infty)$ the
class $F\subset \{f:S\to [0,\infty]~|~\int f(x) dx=1, ~|f(x) - f(y)|\leq
\kappa\|x-y\|_2 ~\forall x,y\in S\}$. Alternatively, if $X^1$ is integrable
and there exist $\alpha,\beta\in \reals$ such that $f(x) \leq \exp(\alpha +
\beta\|x\|_\infty)$ for all $f\in F$, then again the expectation assumption in the corollary is
satisfied.

We next turn the attention to LS regression.  Suppose that we are given the
random design model
\[
Y^j = f^0(X^j) + Z^j, ~~j=1, 2, \dots ,
\]
where the iid random vectors $X^1, X^2, \dots $ take values in the closed set
$S\subset\reals^d$, the iid zero-mean and finite-variance random variables $Z^1, Z^2, \dots $ are
also independent of $X^1, X^2, \dots $, and $f^0:S\to \reals$ is an unknown
function to be estimated based on observations of $(X^1,Y^1)$. Let
\[
L^2_{P}(f,g) = \int \big(f(x) - g(x)\big)^2 dP(x),
\]
where $P$ is the distribution of $X^1$. Consistency in the aw-distance is stated next; see \cite{GeerWegkamp.96} for consistency in the empirical $L^2$ sense.

\begin{corollary}{\rm (consistency in LS regression).}\label{cor:conLS}
Suppose that $\{\epsilon^n\geq 0, n\in\nats\}\to 0$ and $F$ is a closed subset of $\uscfcns(S)$ equi-usc at every $x\in S$. For the random design model above and
\begin{equation*} \hat f^n \in
\epsilon^n\mbox{-}\nargmin_{f\in F} \frac{1}{n}\sum_{j=1}^n \big(Y^j - f(X^j)\big)^2,
\end{equation*}
we have, almost surely, that every cluster point $f^\star$ of $\{\hat f^n, n\in\nats\}$ satisfies
\[
f^\star \in \nargmin_{f\in F} L_P^2(f,f^0).
\]
If $\inf_{f\in F} \Ex[(Y^1 - f(X^1))^2]<\infty$, which occurs in particular when $f^0 \in F$, then $\{\hat f^n, n\in\nats\}$ has at least one cluster point.

When $f^0\in F$, we also have that $f^\star(x) = f^0(x)$ for $P$-a.e. $x\in S$.
\end{corollary}

We next turn to consistency in the presence of sieves, i.e., the class of
functions $F^n$ varies with $n$. The importance of sieves is well-documented
and prior studies include
\cite{Devore.77a,Devore.77b,Grenander.81,GemanHwang.82,DechevskyPenev.97,Chen.07};
see also \cite[Thms. 8.4 and 8.12]{vanderVaart.11}.

\begin{theorem}\label{thm:con_approx}{\rm (consistency; sieves).} Suppose that $X^1, X^2, \dots$ are iid random vectors with values in
$S^0\subset\reals^{d_0}$, $F$ is a closed subset of $\uscfcns(S)$,
$F^n\subset F$, $\psi:S^0\times F\to \Reals$ is a locally
inf-integrable random lsc function, $\pi^n:F\to [0,\infty)$ satisfies
$\pi^n(f^n) \to 0$ for every convergent sequence $\{f^n \in F, n\in \nats\}$,
and $\delta>0$. If $\{\epsilon^n\geq 0, n\in\nats\}\to 0$, then
\begin{align*}
&\nOutLim \Big(\epsilon^n\mbox{-}\nargmin_{f\in F^n_\delta} \frac{1}{n}\sum_{j=1}^n \psi(X^j,f) + \pi^n(f)\Big)\\
&~~~~~~~~~~~~~~~~~~~~ \subset \big\{f\in F^\infty_\delta ~\big|~ \Ex[\psi(X^1,f)] \leq \ninf_{g\in \nLim F^n} \Ex[\psi(X^1,g)]\big\}  \mbox{ a.s.},
\end{align*}
where $F^n_\delta = \{f \in F~|\ninf_{g\in F^n} \setd(f,g)\leq \delta\}$ and
$F^\infty_\delta$ is defined similarly with $F^n$ replaced by $\nLim F^n$. In
particular, if $\nLim F^n = F$, then the right-hand side of the inclusion
equals $\nargmin_{f\in F} \Ex[\psi(X^1,f)]$.
\end{theorem}

The assumptions of the theorem are nearly identical to those of Theorem
\ref{thm:con}. The main difference is that consistency is ensured for
estimators that are near-minimizers of a slightly {\it relaxed} problem over
the class $F^n_\delta$ and not over $F^n$. This relaxation is potentially
beneficial from a computationally point of view (see Section 5.1).

Theorem \ref{thm:con_approx} guarantees that estimators selected from such
relaxed classes will be consistent in some sense. Specifically, every cluster
point of the estimators is at least as ``good'' as $\ninf_{g\in \nLim F^n}
\Ex[\psi(X^1,g)]$ and is also in $F^\infty_\delta$. If $F^n$ eventually
``fills'' $F$, consistency takes place in the usual sense.

To illustrate one application area, we specialize the theorem for ML
estimation of densities, while retaining some of its notation.

\begin{corollary}{\rm (consistency in ML estimation; sieves).}\label{cor:conML_approx}
Suppose that $X^1, X^2$, $\dots$ are iid random vectors, each distributed according to a density
$f^0:S\to[0,\infty]$, $F$ is a closed subset of $\uscfcns(S)$ consisting of
densities, $F^n\subset F$, and for every $f\in F$ there exists $\rho>0$ such
that $\Ex[\nsup_{g\in F} \{\log g(X^1)~|~\setd(f,g)\leq \rho\}] <\infty$. If
$\delta>0$, $\{\epsilon^n\geq 0, n\in\nats\}\to 0$, $f^0\in \nLim F^n$, and
\begin{equation*} \hat f^n \in
\epsilon^n\mbox{-}\nargmin_{f\in F_\delta^n} -\frac{1}{n}\sum_{j=1}^n \log  f(X^j),
\end{equation*}
then, almost surely, $\{\hat f^n, n\in\nats\}$ has at least one cluster point
and every such point $f^\star$ satisfies
\[
K(f^0; f^\star) = 0 \mbox{ and } f^\star \in F_\delta^\infty.
\]
Thus, $f^\star(x) = f^0(x)$ for Lebesgue-a.e. $x\in S$.
\end{corollary}

\subsection{Plug-In Estimators}

Among the many plug-in estimators that can be constructed from density
estimators, those of modes, near-modes, height of modes, and high-likelihood
events are especially accessible within our framework because strong
consistency is {\it automatically} inherited from that of the density
estimator. Similarly, plug-in estimators of ``peaks'' of regression functions and level-sets of classifiers
will also be consistent. Maxima and maximizers of regression functions are
important, especially in engineering design where ``surrogate models'' are
built using regression and that are subsequently maximized to find an optimal
design or decision.

We recall that $\epsilon\mbox{-}\nargmax_{x\in S} f(x) = \{y\in S~|~f(y) \geq
\nsup_{x\in S} f(x) - \epsilon\}$ for $\epsilon\geq 0$ and $f:S\to
\Reals$. Thus, $f(x^\star) = \infty$ when $x^\star \in
\epsilon\mbox{-}\nargmax_{x\in S} f(x)$ and $\nsup_{x\in S} f(x) = \infty$.
If $f$ is a density, then $\nargmax_{x\in S}
f(x)$ is the set of {\it modes} of $f$, $\delta\mbox{-}\nargmax_{x\in S}$
$f(x)$ is a set of {\it near-modes}, and $\nlev_{\geq\alpha} f$ is a set of {\it
high-likelihood events}. We stress that modes are defined here as {\it
global} maximizers of densities. Extension to a more inclusive definition is
possible but omitted.

\begin{theorem}\label{thm:plug}{\rm (plug-in estimators of modes and related quantities).} Suppose that
estimators $\hat f^n\to f^0$ almost surely, with estimates being functions in
$\uscfcns(S)$. If $\{\delta^n\geq 0, ~n\in\nats\}\to \delta$ and
$\{\alpha^n\in \Reals, ~n\in\nats\}\to \alpha$, then the plug-in
estimators
\[
\hat m^n \in\delta^n\mbox{-}\nargmax_{x\in S} \hat f^n(x) ~ \mbox{ and } ~ \hat l^n \in \nlev_{\geq \alpha^n} \hat f^n
\]
are consistent in the sense that  almost surely $\delta\mbox{-}\nargmax_{x\in
S} f^0(x)$ and $\nlev_{\geq \alpha} f^0$ contain every cluster point of
$\{\hat m^n, n\in\nats\}$ and $\{\hat l^n, n\in\nats\}$, respectively.

Moreover, if there is a compact $B\subset S$ such that for all $n$
$\nargmax_{x\in S} \hat f^n(x) \cap B$ $\neq \emptyset$ almost surely, then the
plug-in estimator
\[
\hat  h^n = \nsup_{x\in S} \hat f^n(x) \to \nsup_{x\in S} f^0(x) \mbox{ almost
surely}.
\]
\end{theorem}

The theorem provides foundations for a rich class of {\it constrained}
estimators for modes, near-modes, height of modes, and high-likelihood events
and similar quantities for regression functions and classifiers. We observe that the theorem
holds even if $f^0$ fails to have a unique maximizer. Convergence of
densities in the sense of $L^1$, $L^2$, Hellinger, and Kullback-Leibler as
well as pointwise convergence fails to ensure convergence of modes and
related quantities without additional assumptions.

\section{Closed Classes}

The central technical challenge associated with applying our existence and
consistency theorems is often to establish that the class of functions under
consideration is a closed subset of $\uscfcns(S)$. The analysis is
significantly simplified by the fact that any intersection of closed sets is
also closed. Thus, it suffices to examine each {\it individual} requirement
of a class separately.

It is well known that the limit of a hypo-converging sequence of concave
functions must also be concave and thus the class of concave functions is
closed \cite[Prop. 4.15]{VaAn}. In this section, we provide numerous results for other classes. We note that $S$ is necessarily convex when $f\in\uscfcns(S)$ is convex, concave, or log-concave.

\begin{proposition}\label{prop:closedness_convex}{\rm (convexity and (log-)concavity).}
 For $\{f, f^n, n\in\nats\}$ $\subset$ $\uscfcns(S)$ and $f^n\to f$, we have:
\begin{enumerate}[(i)]

\item If $\{f^n, n\in\nats\}$ are concave, then $f$ is concave. Moreover, if the functions are
    finite-valued, $\kappa\geq 0$, and $\|v\|_2 \leq \kappa$ for
    every subgradient $v\in
    \partial f^n(x)$ and $x\in S$, then $\|v\|_2 \leq \kappa$ for
    every $v\in \partial f(x)$ and $x\in S$.

\item If $\{f^n\geq 0, n\in\nats\}$ are log-concave, then $f$ is
    log-concave.

\item If $\{f^n, n\in\nats\}$ are convex and $\nt S$ is nonempty, then $f$
    is  convex.

 \end{enumerate}
\end{proposition}

The additional assumption about $\nt S$ being nonempty for the convex case is caused by the fact that the aw-distance is inherently tied to hypographs, which makes the treatment of convex functions slightly more delicate than that of concave functions.

Transformations of convex and concave functions beyond the log-concave case lead to the rich
class of s-concave densities; see for example
\cite{SereginWellner.10,KnkM10}.

 \begin{proposition}\label{prop:closedness_trans}{\rm (monotone transformations).}
 For a continuous nondecreasing function $h_0:\reals\to \Reals$, let $h:\Reals \to \Reals$ have $h(y) = h_0(y)$ if $y\in\reals$, $h(-\infty)=\inf_{\bar y\in\reals} h_0(\bar y)$, and $h(\infty)=\sup_{\bar y\in\reals} h_0(\bar y)$. Then, for $\{g^n:S\to \Reals, n\in\nats\}$, with $h\comp g^n\in \uscfcns(S)\to f\in\uscfcns(S)$, the following hold:
 \begin{enumerate}[(i)]

\item If $\{g^n, ~n\in\nats\}$ are concave, then $f=h\comp g$ for some
    concave $g:S\to\Reals$.

\item If $\{g^n, ~n\in\nats\}$ are convex and $\nt S$ is nonempty, then
    $f=h\comp g$ for some convex $g:S\to\Reals$.
 \end{enumerate}
\end{proposition}
Since $h\comp g$ with $h$ nonincreasing and $g$ convex can be written as $\tilde h\comp \tilde g$ with $\tilde h$ nondecreasing and $\tilde g$ concave, the proposition also addresses nonincreasing functions and in fact all {\it
s-concave} functions. This ensures closedness for classes of functions under such shape
restrictions.

 \begin{proposition}\label{prop:closedness_mono}{\rm (monotonicity).}
 For $\{f, f^n, n\in\nats\}\subset \uscfcns(S)$ and $f^n\to f$, we have:
\begin{enumerate}[(i)]

\item If $f^n$ is nondecreasing in the sense that $f^n(x)\leq f^n(y)$ for
    $x\in S$, $y\in \nt S$, with $x\leq y$, then
    $f$ is also nondecreasing in the same sense.

    If $S$ is a box\footnote{A box in $\reals^d$ is of the form $S=[\alpha_1, \beta_2] \times \dots [\alpha_d,
    \beta_d]$, with $-\infty \leq \alpha_i < \beta_i \leq \infty$, where in
    the case of $\alpha_i = -\infty$ and $\beta_i = \infty$ the closed
    intervals are replaced by (half)open intervals. Its dimension is therefore $d$.}, then $\nt S$ can be
    replaced by $S$.

\item If $f^n$ is nonincreasing in the sense that $f^n(x)\geq f^n(y)$ for
    $x\in \nt S$, $y\in S$, with $x\leq y$, then $f$ is also nonincreasing
    in the same sense.

    If $S$ is a box, then $\nt S$ can be replaced by $S$.

 \end{enumerate}
\end{proposition}

The limit of a hypo-converging sequence of nondecreasing functions is not
necessarily nondecreasing for arbitrary $S$. Consider $S = \{(x_1,x_2)\in
\reals^2~|~ x_1 = x_2, 0\leq x_1, x_2\leq 1\} \cup \{(2,0)\}$, $f(x) = f^n(x)
= 0$ if $x=(2,0)$, and $f(x) = 1$ and $f^n(x) = \min\{1, n(x_1+x_2)\}$
otherwise. Clearly, $x = (0,0) \leq y = (2,0)$, but $f(x) = 1 > f(y) = 0$.
Meanwhile, $f^n(x) = f^n(y) = 0$ for all $n$ at these two points and it is
nondecreasing elsewhere too. Still, $f^n\to f$.

We recall that $f:S\to \Reals$ is Lipschitz continuous with modulus
$\kappa$ when $|f(x) - f(y)| \leq \kappa\|x-y\|$ for all $x,y\in S$.

 \begin{proposition}\label{prop:closedness_lip}{\rm (Lipschitz continuity).}
Suppose that $\{f, f^n, n\in\nats\}\subset \uscfcns(S)$, $f^n\to f$, and
$\{f^n, n\in\nats\}$ are Lipschitz continuous with common modulus $\kappa$.
Then, $f$ is also Lipschitz continuous with modulus $\kappa$.
\end{proposition}

 \begin{proposition}\label{prop:closedness_bd}{\rm (pointwise bounds).}
Suppose that $g:S\to \Reals$, $h\in \uscfcns$ $(S)$, $\{f, f^n,
n\in\nats\}\subset \uscfcns(S)$, and $f^n\to f$. If $g(x) \leq f^n(x) \leq h(x)$ for all $n\in\nats$ and $x\in S$, then $g(x) \leq f(x) \leq h(x)$ for all $x\in S$.
\end{proposition}

A function $f:S\to \Reals$ is in the class of multivariate totally
positive functions of order two when $f(x)f(y) \leq
f(\min\{x,y\})f(\max\{x,y\})$ for all $x,y\in S$; see for example
\cite{FallatEtAl.17}. The min and max are taken componentwise.

\begin{proposition}\label{prop:mtp2}{\rm (multivariate total positivity of order two).}
If $\{f^n, n\in\nats\}\subset \uscfcns(S)$ is equi-usc at $\bar x\in S$, the functions
$f^n$ are multivariate totally positive of order two, and $f^n\to f\in \uscfcns(S)$, then
$f$ is multivariate totally positive of order two.
\end{proposition}

Penalty terms and constraints are often defined in terms of sup-functions and integrals. Their (semi)-continuity properties are recorded next.

\begin{proposition}\label{prop:lscsupnorm}{\rm (lsc of sup-norm).}
If $F\subset \uscfcns(S)$ and $g:\Reals\to \Reals$ is lsc\footnote{$g:\Reals\to \Reals$ is lsc if $\nliminf g(y^n) \geq g(y)$ for every $y^n\to y\in \Reals$.}, then $\pi:F\to \Reals$ defined by $\pi(f) = \sup_{x\in S} g(f(x))$ is lsc.
\end{proposition}

In particular, $f\mapsto \nsup_{x\in S} |f(x)|$ is lsc because this corresponds to having $g(y) = |y|$ for $y\in \reals$ and $g(y) = \infty$ for $y = -\infty$ and $\infty$ in the proposition.

\begin{proposition}\label{prop:integral2}{\rm (integral quantities).}
If $\{f^n, n\in\nats\}\subset \uscfcns(S)$ is equi-usc at Lebesgue-a.e. $x\in S$, $f^n\to f\in\uscfcns(S)$, and
for some (measurable) $g:S\to [0,\infty]$, $|f^n(x)|\leq g(x)$ for all $x\in S$
and $n\in\nats$, then
\begin{enumerate}[(i)]

\item $\int f^n(x) dx \to \int f(x) dx$ provided $\int g(x)dx<\infty$;

\item $\int xf^n(x) dx \to \int xf(x) dx$ provided $\int \|x\|g(x)dx
    <\infty$.

\end{enumerate}
\end{proposition}

We end the section with an example of approximating and/or evolving moment information in the definition of a function class.

\begin{proposition}{\rm (moment information.)}\label{prop:moment} Suppose that
$C\subset C^n\subset \reals^d$ are closed, $F^0\subset\uscfcns(S)$ is closed and
equi-usc at every $x\in S$, and there is a function $g:S\to [0,\infty]$ with $\int
\|x\| g(x)dx < \infty$ and $|f(x)|\leq g(x)$ for all $x\in S$ and $f\in
F^0$. Let
\[
F= \Big\{f\in F^0~\Big|~\int xf(x)dx \in C\Big\} \mbox{ and } F^n = \Big\{f\in F^0~\Big|~\int x f(x) dx \in C^n\Big\}.
\]
If $C^n$ set-converges to $C$, then $F^n$ set-converges to $F$.
\end{proposition}

\section{Estimation Algorithm}
For given data $x^1, \dots, x^n\in S^0\subset \reals^{d_0}$, there are no general algorithms available for finding a function in
\begin{equation}\label{eqn:estimatorProblem1}
\epsilon\mbox{-}\nargmin_{f\in F} \frac{1}{n} \sum_{j=1}^n \psi(x^j,f) + \pi(f).
\end{equation}
In this
section, we provide an algorithm for this purpose that combines the need for
approximation of functions in $\uscfcns(S)$ with the use of state-of-the-art
solvers for finite-dimensional optimization.

Suppose that $\pi^\nu$ is an approximation of $\pi$ and $F^\nu$ is an
approximation of $F$ involving only functions that are described by a {\it
finite} number of parameters, i.e., $F^\nu$ is a parametric class. (The sample size $n$ is fixed and we therefore let $\nu\in\nats$
index sequences.) We assume that the statistician finds the
class $F$ appropriate and, for example, believes it balances over- and
underfitting. Consequently, the goal becomes to find a
function in (\ref{eqn:estimatorProblem1}). The approximation $F^\nu$ is
introduced for computational reasons and is
 often selected as
close to $F$ as possible, only limited by the computing resources available.\\

\state Estimation Algorithm.

\begin{description}

  \item[Step 0.]  Set $\nu = 1$.

  \item[Step 1.]  Find $f^\nu \in \epsilon^\nu\mbox{-}\nargmin_{f\in F^\nu}
      \frac{1}{n} \sum_{j=1}^n \psi(x^j,f) + \pi^\nu(f)$.

  \item[Step 2.]  Replace $\nu$ by $\nu +1$ and go to  Step 1.

\end{description}

This seemingly simple algorithm captures a large variety of situations. It
constructs a sequence of functions that approximate those in
(\ref{eqn:estimatorProblem1}) by allowing a tolerance $\epsilon^\nu$ that may
be larger than $\epsilon$ and by resorting to approximations $F^\nu$ and
$\pi^\nu$ of the actual quantities $F$ and $\pi$. The difficulty in carrying
out Step 1 depends on many factors, but since $F^\nu$ consists only of
functions described by a finite number of parameters it reduces to
finite-dimensional optimization for which there are a large number of solvers
available. Section 5.2 shows that we often end up with {\it convex problems.}

The algorithm permits the strategy of initially considering coarse
approximations in Step 1 with subsequent refinement. Since iteration number
$\nu$ has $f^{\nu-1}$ available for warm-staring the computations of $f^\nu$,
the amount of computational work required by a solver in Step 1 is often low.
In essence, the algorithm can make much progress towards
(\ref{eqn:estimatorProblem1}) using relatively coarse approximations.

\begin{theorem}\label{thm:algo}{\rm (convergence of algorithm).} Suppose that
$x^1, \dots, x^n\in \reals^{d_0}$, $F^\nu,F\subset F^0\subset \uscfcns(S)$ are
closed, $\psi(x^j,\cdot):F^0\to (-\infty,\infty]$ is continuous for all $j$,
and $\pi,\pi^\nu:F^0\to \reals$ satisfy $\pi^\nu(g^\nu) \to \pi(g)$ whenever
$g^\nu\in F^0\to g$. Moreover, let $\{\epsilon^\nu\geq 0, \nu\in\nats\}\to
\epsilon^\infty$, $\Lim F^\nu = F$, and $\{f^\nu, \nu\in\nats\}$ be  generated by
the Estimation Algorithm.
\begin{enumerate}[(i)]

\item If $\epsilon^\infty \leq \epsilon$, then
    (\ref{eqn:estimatorProblem1}) contains every cluster point of $\{f^\nu,
    \nu\in\nats\}$.

\item If $\epsilon^\infty < \epsilon$, $F^\nu\subset F$, $F^0$ is bounded, and there exists
    $g\in F$ such that $\psi(x^j,g)<\infty$ for all $j$, then
    (\ref{eqn:estimatorProblem1}) contains $f^{\bar\nu}$ for some finite
    $\bar \nu$.

\end{enumerate}
\end{theorem}
When $\epsilon>0$, item (ii) of the theorem establishes that we obtain an
estimate in a {\it finite} number of iterations of the Estimation Algorithm
as long as $F^\nu$ approximates $F$ from the ``inside.'' Although not the
only possibility, such inner approximations are the primary forms as seen in
Section 5.1.

The main technical and practical challenge associated with the Estimation
Algorithm is the construction of  a parametric class $F^\nu$ that
set-converges to $F$. Since $F$ can be a rich class of usc functions,
standard approaches (see for example
\cite{Meyer.11,Myer12:splines,Myer12:shape}) may fail and we leverage instead
a tailored approximation theory for $\uscfcns(S)$.

\subsection{Parametric Class of Epi-Splines}

Epi-splines is a parametric class that is dense in $\uscfcns(S)$ after a sign
change and furnish the building blocks for constructing a parametric class
$F^\nu$ that approximates $F$. In essence, an epi-spline on
$S\subset\reals^d$ is a piecewise polynomial function that is defined in
terms of a partition of $S$ consisting of $N$ disjoint open subsets that is
dense in $S$. On each such subset, the epi-spline is a polynomial function.
Outside these subsets, the epi-spline is defined by the lower limit of
function values making epi-splines lsc; see \cite{RoysetWets.15b,Royset.16,RoysetWets.15c}. Although approximation theory for epi-splines exists for noncompact
$S$, arbitrary partitions, and higher-order polynomials, we here develop the
possibilities in the statistical setting for a compact polyhedral $S\subset
\reals^d$, simplicial complex partitions, and first-degree polynomials.

We denote by $\cl A$ the closure of a set $A\subset \reals^d$. A collection
 $\cR = \{R_k\}_{k=1}^N$ of open subsets of $S$ is a {\it simplicial complex partition} of $S$ if $\cl R_1$, \dots, $\cl R_{N}$ are
 simplexes\footnote{A {\it simplex} in $\reals^d$ is the convex hull of $d+1$ points $x^0, x^1,
\dots, x^d\in \reals^d$, with $x^1-x^0$, $x^2-x^0$, \dots, $x^d-x^0$ linearly
independent.}, $\cup_{k=1}^{N} \cl R_k = S$,  and $R_k \cap R_l = \emptyset,
k\neq l$. Suppose that $\{\cR^\nu = (R_1^\nu, \dots, R_{N^\nu}^\nu),
\nu\in\nats\}$ is a collection of simplicial complex partition of $S$ with
mesh size $\nmax_{k=1, \dots, N^\nu} \nsup_{x,y\in R_k^\nu} \|x-y\|\to 0$
as $\nu\to \infty$.

A {\it first-order epi-spline} $s$ on a simplicial complex partition $\cR =
\{R_k\}_{k=1}^N$ is a real-valued function that on each $R_k$ is affine and
that satisfies $\nliminf s(x^\nu) = s(x)$ for all $x^\nu\to x$. Let
$\espl(\cR)$ be the collection of all such epi-splines. We deduce from
\cite{RoysetWets.15b,Royset.16} that
\[
\bigcup_{\nu\in\nats} \big\{f:S\to \reals~\big|~f = -s, ~s\in \espl(\cR^\nu)\big\} \mbox{ is dense in } \big(\uscfcns(S),\setd\big).
\]

In the context of the Estimation Algorithm and Theorem \ref{thm:algo}, this
fact underpins several approaches to constructing a parametric class $F^\nu$
that set-converges to $F$. For example, suppose that $F$ is solid\footnote{A
set $A$ is solid if $\cl (\nt A) = A$.}, then $F^\nu = F\cap \espl(\cR^\nu)
\to F$ as can be established by a standard triangular array argument.
One particular class of functions that always will be solid is $F^n_\delta$
in Theorem \ref{thm:con_approx} provided that it is a subset of a convex
$F^0$. For example, $F^0$ can be taken to be $\{f\in \uscfcns(S)~|~f(x)\geq
\alpha ~\forall x\in S\}$, which is convex, so this is no real limitation.
Consequently, the relaxation of $F^n$ to $F_\delta^n$ in Theorem
\ref{thm:con_approx} not only facilitates consistency of an estimator, it
also supports the development of computational methods.

\subsection{Examples of Formulations}

If $F^\nu$ is defined in terms of first-order epi-splines on a partition of
$S\subset\reals^d$ consisting of $N^\nu$ open sets, then each function in
$F^\nu$ is characterized by $N^\nu(d+1)$ parameters. Consequently, Step 1 of
the Estimation Algorithm amounts to approximately solving an optimization
problem with $N^\nu(d+1)$ variables. The number of variables is independent
of the sample size $n$. The number of open sets $N^\nu$ would usually grow
with $d$, but when the growth is slow the number of variables is manageable
for modern optimization solvers even for moderately large $d$.

Among the numerous formulations of the problem in Step
1 of the Estimation Algorithm, we illustrate one based on first-order
epi-splines with a simplicial complex partition, which is also used in
Section 2.3. Suppose that $c_k^0, c_k^1, \dots, c_k^d\in \reals^d$ are the vertexes of the
$k$th simplex of a simplicial complex partition of $S\subset\reals^d$ with
$N$ simplexes. A first-order epi-spline is then fully defined by its height
at these vertexes. Let $h_k^i\in \reals$ be the height at $c_k^i$, $i=0, 1,
\dots, d$, $k=1, \dots, N$. These $N(d+1)$ variables are to be optimized.
(Optimization over such ``tent poles'' is familiar in ML estimation over
log-concave densities, but then they are located at the data points and not
according to simplexes as here; see for example
\cite{CuleSamworthStewart.10a}.) We next give specific expressions for
typical objective and constraint functions.

In ML estimation of densities, the loss expressed in terms of the
optimization variables becomes
\[
-\frac{1}{n} \sum_{j=1}^n \log f(x^j) = -\frac{1}{n} \sum_{j=1}^n\log \sum_{i=0}^d \mu^i_{j} h_{k_j}^i,
\]
where $k_j$ is the simplex in which data point $x^j$ is located and the
scalars $\{\mu^i_{j}, i=0, 1, \dots, d, j=1, \dots, n\}$ can be precomputed
by solving $x^j = \sum_{i=0}^d \mu^i_{j} c_{k_j}^i$. The loss is therefore
convex in the optimization variables.

The requirement that functions are nonnegativity is implemented by the
constraints $h_k^i \geq 0$  for all $i=0, 1, \dots, d, ~k=1, \dots, N$,
which define a polyhedral feasible set.

The requirement that functions integrate to one is implemented by
\[
\int f(x) dx = \frac{1}{d+1}\sum_{k=1}^N \alpha_k \sum_{i=0}^d h_{k}^i =1,
\]
where $\alpha_k$ is the hyper-volume of the $k$th simplex.

The requirement that functions should have their argmax covering a given
point $x^\star$ is implemented by the constraints
\[
\sum_{i=0}^d \eta^i h_{k^\star}^i \geq h_k^{i'} \mbox{ for all } i'=0, 1,
\dots, d, ~k=1, \dots, N,
\]
where $k^\star$ is the simplex in which $x^\star$ is located and the scalars
$\{\eta^i, i=0, 1, \dots, d\}$ can be precomputed by solving $x^\star =
\sum_{i=0}^d \eta^i c_{k^\star}^i$. The constraints form a polyhedral
feasible set.

Implementation of continuity, Lipschitz continuity, concavity, and many other
conditions also lead to polyhedral feasible sets. Consequently, ML estimation
of densities on a compact polyhedral set $S\subset\reals^d$ under a large
variety of constraints can be achieved by optimization of a convex function
over a polyhedral feasible sets for which highly efficient solvers are available. A switch to LS regression, would result in
a convex quadratic function to minimize, with many of the constraints
remaining unchanged. In that case, specialized quadratic optimization solvers
apply.\\

\noindent {\bf Acknowledgements.} This material is based upon work supported in part by ONR Science of Autonomy (N0001417WX01210, N000141712372), DARPA (HR0011834187), and NPS CIMS.\\

\bibliographystyle{plain}
\bibliography{refs}

\section{Appendix: Additional Examples}

This section discusses existence of solutions of approximation problems for the already well-understood classes of monotone and of log-concave functions. We give proofs passing through the metric space $(\uscfcns(S),\setd)$ to further illustrate the framework.

\begin{proposition}\label{prop:exist_LS-abscont}{\rm (existence of monotone LS approximation).}
For a box $S\subset\reals^d$, suppose that $F = \{f\in \uscfcns(S)~|~f \mbox{ nondecreasing}\}$ and $P$ is an absolutely continuous distribution on $S\times\reals$. Then,
\[
\nargmin_{f\in F} \int (y-f(x))^2 dP(x,y) \neq\emptyset.
\]
\end{proposition}
\state Proof. By Proposition \ref{prop:closedness_mono}, $F$ is closed. Suppose that $f^n\in F\to f$. Let $D = \{x\in \nt S~|~f \mbox{ is discontinuous at } x\}$. In view of Propositions \ref{prop:equi-usc}(iv) and \ref{prop:equi-usc_ptwise}, $f^n(x) \to f(x)$ for all $x\in \nt S\setminus D$. Thus, $(y-f^n(x))^2 \to (y-f(x))^2$ for all such $x$ and all $y\in \reals$.
By \cite{Lavric.93}, $D$ has Lebesgue measure zero and the same holds for $S\setminus \nt S$. Then, by Fatou's Lemma, $\nliminf \int (y-f^n(x))^2 dP(x,y) \geq \int (y-f(x))^2 dP(x,y)$ and $f\mapsto \int (y-f(x))^2 dP(x,y)$ is lsc on $F$. Its lower level-sets are compact at every finite level (cf. the argument in the proof of Corollary \ref{cor:existence-np}) and the conclusion follows.\eop

The next result is in \cite{DumbgenSamworthSchuhmacher.11}, but we provide a proof with some novel elements: the log-likelihood criterion function is shown to be lsc on the enlarged class of log-concave functions that integrate to values in $[0,1]$.

\begin{proposition}\label{prop:exist_ML-logconcave}{\rm (existence of log-concave ML estimator).}
Suppose that $F=\{f\in \uscfcns(\reals^d)~|~f\geq 0, \mbox{ log-concave}\}$. Then, for any probability distribution $P$ on $\reals^d$,
\[
\nargmin_{f\in F} \Big\{\int -\log f(x) dP(x)~\Big|~\int f(x) dx = 1\Big\} \neq\emptyset
\]
if and only if
\[
\int \|x\|dP(x) < \infty \mbox{ and } P(H)<1 \mbox{ for all hyperplane } H\subset\reals^d.
\]
\end{proposition}
\state Proof. For $f^n\in F\to f$, Proposition \ref{prop:closedness_convex}(ii) establishes that $f$ is log-concave. Moreover, $f^n(x)\to f(x)$ for all $x\in \nt \{x\in \reals^d~|~f(x)>0\}$ and also when $f(x) = 0$ by Propositions \ref{prop:equi-usc_ptwise} and \ref{prop:equi-usc}. The subset of $\reals^d$ that fails outside both of these cases has Lebesgue measure zero so $f^n(x)\to f(x)$ for Lebesgue-a.e. $x\in\reals^d$. Fatou's Lemma gives that $\nliminf \int f^n(x) dx$ $\geq$ $\int f(x) dx$. Thus, $F_\leq = \{f\in F~|~\int f(x) dx \leq 1\}$  is closed and actually compact because all functions in $F$ are nonnegative.

We show that $\phi(f) = \int -\log f(x)dP(x)$ is lsc as a function on $(F_\leq,\setd)$. Let $f^n\in F_\leq \to f$. We consider two cases: a) $\int f(x) dx = \gamma > 0$. Then, $\gamma^{-1} f$ is a log-concave density and by \cite[Lem. 1]{CuleSamworth.10} there are $\xi_0\in\reals$ and $\xi_1\in (0,\infty)$ such that $f(x) \leq \exp(\xi_0-\xi_1\|x\|)$ for all $x\in\reals^d$. Let $\epsilon = \nsup_{x\in\reals^d} f(x)/4$, which then must be positive, and $\rho \in (2\epsilon, \infty)$ such that $f(x) \leq \epsilon$ for $\|x\|_2\geq \rho/2$. (Here, we adopt the Euclidean norm, with the correspond balls denoted by $\ball_2(x,\delta)$, to simplify a reference to \cite{VaAn}.) Hypo-convergence is locally uniform in the following sense \cite[Thm. 4.10]{VaAn}: there is $\bar n$ such that for $n\geq \bar n$,
\begin{align*}
  \hypo f^n \cap \ball_2(0,\rho) & \subset \hypo f + \ball_2(0,\epsilon)\\
  \hypo f \cap \ball_2(0,\rho) & \subset \hypo f^n + \ball_2(0,\epsilon).
\end{align*}
Take $(x,f^n(x))$ with $\|x\|_2= \rho$. If $f^n(x) >\rho$, then $(x,\rho)\in \hypo f^n \cap \ball_2(0,\rho)$ and there exists $(y,\beta)\in \hypo f$ such that $\|x-y\|_2 \leq \epsilon$ and $|\rho - \beta|\leq \epsilon$. Thus, $f(y) \geq \beta\geq \rho - \epsilon > \epsilon$. However, $f(y) \leq \epsilon$ because $\|y\|_2\geq \rho/2$ and we have reached a contradiction. Thus, $f^n(x) \leq \rho$, $(x,f^n(x))\in \hypo f^n \cap \ball_2(0,\rho)$, and there is $(y,\beta) \in \hypo f$ such that $\|x-y\|_2\leq \epsilon$ and $|f^n(x) - \beta| \leq \epsilon$. This leads to $f^n(x) \leq \beta + \epsilon \leq f(y) + \epsilon \leq 2\epsilon$ for all $n\geq \bar n$. The choice of $\rho$ ensures that $\bar x\in\nargmax_{x\in \reals^d} f(x)$ with $\|\bar x\|_2 \leq \rho/2$ exists. By \eqref{eqn:hypoconv}, there is $x^n\to \bar x$ such that $f^n(x^n) \to f(\bar x) = 4\epsilon$. Thus, for some $n^*\geq \bar n$, $\|x^n\|_2 \leq 3\rho/4$ and $f^n(x^n) \geq 3\epsilon$ for all $n\geq n^*$. Since we also have $f^n(x) \leq 2\epsilon$ for $\|x\|_2 = \rho$, $\nargmax_{x\in \reals^d} f^n(x) \subset \ball_2(0,3\rho/4)$ for all $n\geq n^*$. By \cite[Thm. 7.31]{VaAn}, this implies that $\nsup_{x\in \reals^d} f^n(x) \to \nsup_{x\in \reals^d} f(x)$. Consequently, for sufficiently large $n$, $\int -\log f^n(x) dP(x) \geq \int -\log [2 \nsup_{\bar x\in \reals^d} f(\bar x)] dP(x)>-\infty$, which then furnishes an integrable lower for application of Fatou's lemma: $\nliminf \int -\log f^n(x) dP(x)$ $\geq$ $\int \nliminf [-\log f^n(x)] dP(x)$. Since $\nliminf$ $-\log f^n(x) \geq -\log f(x)$ for all $x\in\reals^d$ by \eqref{eqn:hypoconv}, we conclude that $\phi$ is lsc at points $f\in F_\leq$ with $\int f(x) dx >0$. This fact holds for any $P$.

Next, we consider b) $\int f(x) dx = 0$ and now it becomes essential to limit the scope to $P$ with the stated properties. Let $D = \{x\in\reals^d~|~f(x)>0\}$, which then has Lebesgue measure zero (because $\int f(x) dx = 0$) and $\nt D = \emptyset$. Since $D$ is also convex by the log-concavity of $f$, it lies in an affine subspace of $\reals^d$ of dimension less than $d$, i.e., $D$ is a subset of some hyperplane $H\subset\reals^d$. Consequently, the first term of
\[
\phi(f) = \int_{x\not\in D} -\log f(x) dP(x) + \int_{x\in D} -\log f(x) dP(x)
\]
integrates to $\infty$ in view of the assumption on $P$. The convention $\infty-\alpha = \infty$ for any $\alpha \in \Reals$ implies that $\phi(f) = \infty$ regardless of the value of the second term. It remains to show that $\phi(f^n)\to \infty$ when $f^n\in F_{\leq}\to f$.
Since $\phi(f^n) = \infty$ when $\int f^n(x) dx = 0$ as just argued, we assume without loss of generality that $\int f^n(x) dx>0$ for all $n$. In fact, those integrals can be assumed to be one because, with $\int f^n(x) dx = \gamma^n$, $\int -\log f^n(x) dP(x) = -\log \gamma^n + \int -\log (f^n(x)/\gamma^n) dP(x)\to \infty$ when the last term tends to $\infty$.

Each $s^n=\nsup_{x\in\reals^d} f^n(x)$, $n\in\nats$, is finite (cf. \cite[Lem. 1]{CuleSamworth.10}), but the sequence could be unbounded. If $\nsup_{n\in\nats} s^n$ is also finite, then
\begin{align*}
\phi(f^n) = &\int_{f(x)=0, f^n(x)\leq 1} -\log f^n(x) dP(x) + \int_{f(x)>0, f^n(x)\leq 1} -\log f^n(x) dP(x)\\
&~~~~ + \int_{f^n(x)>1} -\log f^n(x) dP(x) \to \infty;
\end{align*}
the first term tends to $\infty$ because $f^n(x) \to 0$ when $f(x) = 0$ by Proposition \ref{prop:equi-usc}(i) and the last term is bounded from below uniformly in $n$. Hence, suppose that $s^n\to\infty$. For $\eta>0$, $\tau^n = \log s^n$, and $\sigma^n = \exp(-\eta \tau^n)$,
\begin{align*}
\int -\log f^n(x) dP(x) & \geq \eta \tau^n P(\reals^d\setminus \nlev_{\geq \sigma^n} f^n) - \tau^n P(\nlev_{\geq \sigma^n} f^n)\\
& = (\eta+1)\tau^n\Big(\frac{\eta}{\eta+1} - P(\nlev_{\geq \sigma^n} f^n) \Big).
\end{align*}
By \cite[Lem. 4.1]{DumbgenSamworthSchuhmacher.11}, the Lebesgue measure of  $\nlev_{\geq \sigma^n} f^n$ is no greater than
\[
(1+\eta)^d (\tau^n)^d \exp(-\tau^n)/ \int_0^{(1+\eta)\tau^n} t^d \exp(-t) dt \to 0
\]
as $s^n$ (and $\tau^n$) tends to infinity for any given $\eta>0$.  Moreover, \cite[Lem. 2.1]{DumbgenSamworthSchuhmacher.11} establishes that $P(\nlev_{\geq \sigma^n} f^n)<\eta/(\eta+1)$ when the Lebesgue measure of  $\nlev_{\geq \sigma^n} f^n$ is sufficiently low and $\eta$ sufficiently high. (This fact relies critically on the assumption on $P$.) Thus, $\int -\log f^n(x) dP(x)\to \infty$ when $s^n\to \infty$ and $\phi$ is lsc (in fact continuous) at $f$ when $\int f(x) dx = 0$.

In summary, we have shown that $\phi$ is lsc on the compact set $F_\leq$. Thus, there exists  $f^\star \in \nargmin_{f\in F_\leq} \phi(f)$. Trivially, there is $f\in F_\leq$ with finite $\phi(f)$, which implies that $\phi(f^\star)<\infty$ and, as argued above, $\int f^\star (x)dx=\gamma>0$. Since $\phi(f^\star/\gamma)\leq \phi(f^\star)$, $f^\star/\gamma \in \nargmin_{f\in F} \{\phi(f)~|~\int f(x) dx=1\}$.

For the necessity of the conditions on $P$ we refer to \cite{DumbgenSamworthSchuhmacher.11}.\eop

\section{Appendix: Intermediate Results and Proofs}

This section includes proofs of all the results in the paper.\\

\state Proof of Proposition \ref{prop:example}. By Proposition \ref{prop:equi-usc}(v), $F$ is equi-usc at all $x\in \reals^d$. Proposition \ref{prop:integral2}(i)
ensures that the integral constraint is closed. Theorem \ref{thm:plug} as well as  Propositions \ref{prop:closedness_lip} and \ref{prop:closedness_bd}
establish that the other constraints are closed too. Consequently, $F$ is compact. Corollary \ref{cor:exist_ML} applies and confirms
(i). Corollary \ref{cor:conML} and the discussion immediately after it establish
(ii). When $f^0\in F$, then every cluster point of $\{\hat f^n, n\in\nats\}$
must deviate from $f^0$ at most on set of Lebesgue measure zero. For
Lipschitz continuous functions this means that the functions must be
identical and (iii) holds.\eop

\state Proof of Proposition \ref{prop:equi-usc}. When $f^n\to f$, it suffices by \cite[Thm. 7.10]{VaAn} to establish that $f^n(\bar x) \to f(\bar x)$. In view of \eqref{eqn:hypoconv}, (i) is trivial. Items (ii,iii) follow by \cite[Thm. 7.17]{VaAn}. For (iv), we only prove the nondecreasing case as a nearly identical argument establishes the conclusion for nonincreasing functions. Let $\epsilon>0$. Since $\bar x\in \nt S$ and $f$ is continuous at $\bar x$, there exist $\bar y \in S$, with $\bar y_i < \bar x_i$ for $i=1, \dots, d$, and $f(\bar y) \geq f(\bar x) - \epsilon$. Moreover, for some $x^n\in S\to \bar y$, $f^n(x^n) \to f(\bar y)$ by \eqref{eqn:hypoconv}. Since $x^n \leq \bar x$ for sufficiently large $n$, $\nliminf f^n(\bar x) \geq \nliminf f^n(x^n) = f(\bar y) \geq f(\bar x) - \epsilon$. Since $\epsilon$ is arbitrary, the conclusion follows because $\nlimsup f^n(\bar x) \leq f(\bar x)$ already by \eqref{eqn:hypoconv}. For (v), consider the definition of equi-usc. The Lipschitz condition ensures that there is $\delta\in (0,\infty)$ with $f^n(x) \leq f^n(\bar x) + \kappa \|x-\bar x\|$ for all $n\in\nats$ and $x\in \ball(\bar x, \delta)$. Let $\epsilon>0$. If $\kappa = 0$, then set $\delta' = \delta$. Otherwise, set $\delta' = \min\{\epsilon/\kappa, \delta\}$. In either case,
$\nsup_{x\in \ball(\bar x,\delta')} f^n(x) \leq f^n(\bar x) + \kappa\delta' \leq f^n(\bar x) + \epsilon$.\eop

\state Proof of Proposition \ref{prop:integral1}. In view of Proposition \ref{prop:equi-usc_ptwise}, the result follows directly from applications of the Dominated Convergence
Theorem.\eop

\state Proof of Theorem \ref{thm:existence}. A trivial generalization of Fatou's Lemma shows that $\int \psi(x,\cdot) dP^0(x)$ is lsc on $F$ (see for example \cite[Appendix]{DongWets.00}). Since for all $f\in F$, $\int \psi(x,f) dP^0(x)$ and $\pi(f)$ exceed $-\infty$, $\int \psi(x,\cdot) dP^0(x) + \pi$ is lsc on $F$. All lsc
functions defined on a compact set attain their infima.\eop

\state Proof of Corollary \ref{cor:existence}. The function $n^{-1}
\sum_{j=1}^n \psi(x^j,\cdot) + \pi$ is lsc on $F$ because each term in the sum involves a lsc function that is never $-\infty$. All lsc
functions defined on a compact set attain their infima.\eop

\state Proof of Corollary \ref{cor:existence-svm}. By Propositions \ref{prop:closedness_convex}(i) and \ref{prop:closedness_bd} as well as Theorem \ref{thm:plug}, $F$ is closed. It is also bounded; see the remark after Corollary \ref{cor:existence}. Since $g>-\infty$ and $x^j\in \nt S$, it is also equi-usc at $x^j$, $j=1, \dots, n$, by Proposition \ref{prop:equi-usc}(ii). Thus, in view of Proposition \ref{prop:equi-usc_ptwise} $f\mapsto \max\{0, 1- y^j f(x^j)\}$ is continuous on $F$ for all $j$ and the conclusion follows by Corollary \ref{cor:existence}.\eop

\state Proof of Corollary \ref{cor:existence-np}. Let $\phi(f) = n^{-1} \sum_{j=1}^n \max\{0, 1- f(x^j)\}$, $f\in F$, and $f^n\in F\to f$. By \eqref{eqn:hypoconv}, $\nlimsup f^n(x^j) \leq f(x^j)$ and $\nliminf (\max\{0, 1-f^n(x^j)\}) \geq \max\{0, 1 - f(x^j)\}$ for all $j=1, \dots, n$, which implies that $\phi$ is lsc on $F$. Since $Z^i$ is open, $\nliminf (\nsup_{x\in Z^i} f^n(x)) \geq \nsup_{x\in Z^i} f(x)$  by \cite[Prop. 7.29]{VaAn}. Thus, $F^0 = \{f\in F~|~\nsup_{x\in Z^i} f(x) \leq 0, i=1, \dots, m\}$ is closed. For $g\in F$, suppose that $\{f^n\in F, n\in\nats\}$ is such that $\setd(f^n,g)\to \infty$. Then, $\hypo f^n$ set-converges to $\emptyset$, which implies $f^n(x^j)\to -\infty$ for all $j$ and $\phi(f^n)\to \infty$. Consequently, $\{f\in F~|~\phi(f) \leq \alpha\}$ is bounded for $\alpha\in\reals$. Since it is also closed by virtue of $\phi$ being lsc, these level-sets are actually compact. A lsc function with compact level-set attains it infimum.\eop

\state Proof of Corollary \ref{cor:exist_ML}. Let $f^n\in F\to f$. By \eqref{eqn:hypoconv}, $\nlimsup f^n(x^j)\leq f(x^j)$
for all $j$ so $\nliminf -\log f^n(x^j)\geq
-\log f(x^j)$. Thus, $f\mapsto -n^{-1} \sum_{j=1}^n \log f(x^j)$ is lsc on $F$. The conclusion then follows by Corollary \ref{cor:existence}.\eop

\state Proof of Corollary \ref{cor:exist_LS}. In view of Proposition \ref{prop:equi-usc_ptwise}, $f\mapsto  \sum_{j=1}^n (y^j - f(x^j))^2$ is continuous on $F$. An argument similar to the one in the proof of Corollary \ref{cor:existence-np} yields that this function has compact level-sets.\eop

The proof of Theorem \ref{thm:con} relies on an lsc-LLN, essentially in
\cite{ArtsteinWets.96,KorfWets.01}, that ensures almost sure epi-convergence
of empirical processes indexed on a polish space. For completeness, we
include the statement as well as a new proof, which is simpler than that in
\cite{ArtsteinWets.96}. It follows the arguments in \cite{KorfWets.01} for
ergodic processes, but takes advantage of the present iid setting. The
statement is made slightly more general than needed without complication.

\begin{proposition}\label{prop:lln} {\rm (lsc-LLN).} Suppose that $(Y,d_Y)$ is a complete separable (polish) metric space, $(\Xi,\cA,P)$ is a
complete\footnote{In view of \cite{KorfWets.01}, the result (but not our
proof) holds without completeness.} probability space, and $\psi:\Xi\times
Y\to\Reals$ is a locally inf-integrable random lsc
function\footnote{The definitions of Section 3.2 carry over to the more
general context here.}. If $\bfxi^1, \bfxi^2, \dots$ is a sequence of iid
random elements that take values in $\Xi$ with distribution $P$, then almost
surely
\begin{equation*}
\frac{1}{n} \sum_{j=1}^n \psi(\bfxi^j,\cdot) \mbox{ epi-converges } \Ex[\psi(\bfxi^1,\cdot)],
\end{equation*}
which is equivalent to having for all $y\in Y$,
\begin{align*}
  \forall y^n\to y, &\nliminf \frac{1}{n} \sum_{j=1}^n \psi(\bfxi^j,y^n) \geq \Ex[\psi(\bfxi^1,y)]\\
  \exists y^n\to y, &\nlimsup \frac{1}{n} \sum_{j=1}^n \psi(\bfxi^j,y^n) \leq \Ex[\psi(\bfxi^1,y)].
\end{align*}
\end{proposition}
\state Proof. A slight generalization of Fatou's Lemma (see
\cite[Appendix]{DongWets.00}) ensures that $\Ex[\psi(\bfxi^1,\cdot)]$ is lsc.
Let $\bar D\subset Y\times\Reals$ be a countable dense subset of
the epigraph $\epi \Ex[\psi(\bfxi^1,\cdot)]$, with $\epi h = \{(y,y_0)\in
Y\times \reals~|~ h(y)\leq y_0\}$, which may be empty. Moreover, let
$D\subset Y$ be a countable dense subset of $Y$ that contains the projection
of $\bar D$ on $Y$ and $Q_+$ be the nonnegative rational numbers. For $y \in
D$ and $r\in Q_+$, we define $\pi_{y,r}:\Xi\to \Reals$ by setting
\begin{equation*}
  \pi_{y,r}(\xi) = \ninf_{y'\in\ball^o(y,r)} \psi(\xi,y') \mbox{ if } r>0 \mbox{ and } \pi_{y,0}(\xi) = \psi(\xi,y),
\end{equation*}
where $\ball^o(y,r) = \{y'\in Y~|~d_Y(y',y)<r\}$. By Theorem 3.4 in
\cite{KorfWets.01}, every such $\pi_{y,r}$ is an extended real-valued random
variable defined on the probability space $(\Xi,\cA,P)$. Since $\psi$ is
locally inf-integrable, it follows that for every $y\in D$ there is a
closed neighborhood $V_y$ of $y$ and $r_y\in (0,\infty)$ such that
\begin{equation*}
  \ball^o(y,r) \subset V_y \mbox{ and } \Ex[\pi_{y,r}] \geq \int \ninf_{y'\in V_y} \psi(\xi,y') dP(\xi) > -\infty \mbox{ for } r \in [0,r_y].
\end{equation*}
Let $(\Xi^\infty,\cA^\infty,P^\infty)$ be the product space constructed from
$(\Xi,\cA,P)$ in the usual manner. For every $y\in D$ and $r \in [0,r_y] \cap
Q_+$, a standard law of large numbers for extended real-valued random
variables (see for example \cite[Thms. 7.1,7.2]{Durrett.96}) ensures
that
\begin{equation*}
  \frac{1}{n}\sum_{j=1}^n \pi_{y,r}(\xi^j) \to \Ex[\pi_{y,r}] \mbox{ for } P^\infty\mbox{-a.e. } (\xi^1, \xi^2, \dots) \in \Xi^\infty.
\end{equation*}
Since $\{\pi_{y,r}~| ~y\in D, r \in [0,r_y] \cap Q_+\}$ is a countable
collection of random variables, there exists $\Xi_0^\infty\subset \Xi^\infty$
such that $P(\Xi_0^\infty) = 1$ and
\begin{equation*}
  \frac{1}{n}\sum_{j=1}^n \pi_{y,r}(\xi^j) \to \Ex[\pi_{y,r}] \mbox{ for all } (\xi^1, \xi^2, \dots) \in \Xi_0^\infty \mbox{ and } y\in D, r \in [0,r_y] \cap Q_+.
\end{equation*}

We proceed by establishing the liminf and limsup conditions of the theorem.
First, suppose that $y^n\to y$. There exist $\bar n^k\in\nats$, $z^k\in D$,
and $r^k \in [0,r_y] \cap Q_+$, $k\in\nats$, such that $z^k\to y$, $r^k\to
0$,
\begin{equation*}
   \ball^o(z^k,r^k) \supset \ball^o(z^{k+1},r^{k+1}), \mbox{ and } y^n\in \ball^o(z^k,r^k) \mbox{ for } n\geq \bar n^k, k \in \nats.
\end{equation*}
We temporarily fix $k$. Then, for $n\geq \bar n^k$ and $(\xi^1, \xi^2, \dots)
\in \Xi_0^\infty$,
\begin{equation*}
  \frac{1}{n} \sum_{j=1}^n \psi(\xi^j,y^n) \geq \frac{1}{n} \sum_{j=1}^n \ninf_{y'\in\ball^o(z^k,r^k)} \psi(\xi^j,y') = \frac{1}{n} \sum_{j=1}^n \pi_{z^k,r^k}(\xi^j) \to \Ex[\pi_{z^k,r^k}].
\end{equation*}
The nestedness of the balls, implies that $\pi_{z^k,r^k} \leq
\pi_{z^{k+1},r^{k+1}}$ for all $k$. Moreover the lsc of $\psi(\xi,\cdot)$
implies that for all $\xi\in \Xi$, $\pi_{z^k,r^k}(\xi)\to \pi_{y,0}(\xi) =
\psi(\xi,y)$.  Thus, in view of the Monotone Convergence Theorem,
$\Ex[\pi_{z^k,r^k}] \to \Ex[\psi(\bfxi^1,y)]$. We have establish that for
$(\xi^1, \xi^2, \dots) \in \Xi_0^\infty$, $\nliminf$ $n^{-1}$ $\sum_{j=1}^n$
$\psi(\xi^j,y^n) \geq \Ex[\psi(\bfxi^1,y)]$.

Second, for every $y\in Y$, we construct a sequence $y^n\to y$ such that for
$(\xi^1, \xi^2, \dots) \in \Xi_0^\infty$, $\nlimsup n^{-1} \sum_{j=1}^n
\psi(\xi^j,y^n) \leq \Ex[\psi(\bfxi^1,y)]$.

Suppose that $y\in D$. Then, the claim holds because for $(\xi^1, \xi^2,
\dots) \in \Xi_0^\infty$
\begin{equation*}
  \nlimsup \frac{1}{n} \sum_{j=1}^n \psi(\xi^j,y) = \frac{1}{n} \sum_{j=1}^n \pi_{y,0}(\xi^j) \to \Ex[\pi_{y,0}] = \Ex[\psi(\bfxi^1,y)].
\end{equation*}
Fix $(\xi^1, \xi^2, \dots) \in \Xi_0^\infty$ and let $h:Y\to
\Reals$ be the unique lsc functions that has as epigraph the set
$\nOutLim \{\epi n^{-1} \sum_{j=1}^n \psi(\xi^j,\cdot)\}$. Thus, the prior
equality is equivalent to having $h(y) \leq \Ex[\psi(\bfxi^1,y)]$, which then
holds for all $y\in D$. Consequently, $\{(y,\alpha) \in Y\times \reals~|~
h(y) \leq \alpha, y\in D\} \subset \epi  \Ex[\psi(\bfxi^1,\cdot)]$. Since $h$
is lsc and $\epi \Ex[\psi(\bfxi^1,\cdot)]$ is closed, we have after taking
the closure on both sides that $\epi h \subset \epi \Ex[\psi(\bfxi^1,\cdot)]$
and also $h(y) \leq \Ex[\psi(\bfxi^1,y)$ for all $y$. By construction of $h$,
this implies that for all $y$ there exists $y^n\to y$ such that $\nlimsup
n^{-1} \sum_{j=1}^n \psi(\xi^j,y^n) \leq \Ex[\psi(\bfxi^1,y)]$ and the
conclusion holds.\eop

\state Proof of Theorem \ref{thm:con}. If $F$ is empty, the results hold
trivially. Suppose that $F$ is nonempty. By \cite[Prop. 4.45, Thm. 7.58]{VaAn}, $(\uscfcns(S),\setd)$
is a complete separable metric space. By
virtue of being a closed subset, $F$ forms another complete separable metric
space $(F,\setd)$. Let $\phi^n(f) = n^{-1}\sum_{j=1}^n \psi(X^j,f) + \pi^n(f)$ and
$\phi(f) = \Ex[\psi(X^1,f)]$, $f\in F$.
Proposition \ref{prop:lln} applied with this metric space establishes that
$n^{-1}\sum_{j=1}^n \psi(X^j,\cdot)$ epi-converges to $\phi$ a.s. Moreover,
for all $f^n \in F\to f$,
\begin{equation*}
\nliminf \phi^n(f^n)  \geq \nliminf \frac{1}{n} \sum_{j=1}^n \psi(X^j,f^n) \geq \phi(f) \mbox{ a.s.}
\end{equation*}
Also, there exists $f^n\in F\to f$ such that
\begin{align*}
\nlimsup \phi^n(f^n)  &\leq \nlimsup \frac{1}{n} \sum_{j=1}^n \psi(X^j,f^n) + \nlimsup \pi^n(f^n) \leq \phi(f) \mbox{ a.s.}
\end{align*}
We have established that $\phi^n$ epi-converges to $\phi$ a.s. If $\phi$ is
improper, which in this case means that $\phi(f)=\infty$ for all $f\in F$,
then item (i) holds trivially because the right-hand side of the inclusion is
the whole of $F$. If $\phi$ is proper, then $\phi^n$ is also proper and
\cite[Prop. 2.1]{Royset.16} applies, which establishes (i).

The additional assumptions in item (ii) imply that both $\phi$ and $\phi^n$
are proper, and also that $\phi^n$ epi-converges tightly $\phi$ because then
$F$ is compact. Thus, \cite[Thm. 3.8]{RoysetWets.16a} applies and item
(ii) is established.\eop

\state Proof of Corollary \ref{cor:consvm}. We deduce from the proof of Corollary \ref{cor:existence-svm} that $F$ is compact and equi-usc at all $x\in \reals^d$. This implies that for all $(x,y)\in\reals^d\times\{-1,1\}$, $f\mapsto \max\{0, 1-yf(x)\}$ is continuous on $F$. Suppose that $f^n\to f$ and $x^n\to x$. By \eqref{eqn:hypoconv}, $\nlimsup f^n(x^n) \leq f(x)$. Thus, $(x,f)\mapsto f(x)$ is usc on $\reals^d \times F$. From this we conclude that $((x,y),f)\mapsto \max\{0, 1-yf(x)\}$ is measurable and a random lsc function. It is trivally locally inf-integrable by virtue of being nonnegative. Theorem \ref{thm:con}(i) therefore
applies and a cluster point $f^\star$ of $\{\hat f^n, n\in\nats\}$, of which one exits due to compactness of $F$, must satisfy the first conclusion a.s. The second conclusion follows by an application of \cite[Prop. 7.7]{VaAn}.\eop

\state Proof of Corollary \ref{cor:conML}. Since $F$ consists of nonnegative functions, it is bounded and in fact compact since closed. Thus, $\{\hat f^n, n\in\nats\}$
must have at least one cluster point. Next, we show that $\psi:S\times F\to
\Reals$ given by $\psi(x,f) = -\log f(x)$ is a random lsc function.
Suppose that $f^n\in F \to f$ and $x^n\in S\to x$, then $\nlimsup f^n(x^n)
\leq f(x)$ and also $\nliminf -\log f^n(x^n) \geq -\log f(x)$, which implies
that $\psi$ is lsc. Measurability then follows directly from the fact that
lower level-sets of lsc functions are closed. Theorem \ref{thm:con}(i) therefore
applies and a cluster point $f^\star$ of $\{\hat f^n, n\in\nats\}$ must
satisfy a.s.
\[
f^\star \in \nargmin_{f\in F}\Ex[-\log f(X^1)] \subset \nargmin_{f\in F} \Ex[\log f^0(X^1)] -\Ex[\log f(X^1)].
\]
The inclusion holds even if $\Ex[\log f^0(X^1)]$ equals $-\infty$ or
$\infty$. The last conclusion of the theorem follows directly from the
properties of the Kullback-Leibler divergence.\eop

\state Proof of Corollary \ref{cor:conLS}. From the proof of Corollary
\ref{cor:exist_LS}, we deduce that $f\mapsto (y-f(x))^2$ is continuous for any $(x,y)\in S\times
\reals$. Moreover, if $f^n\in F\to f$ and $x^n\in S\to x$, then $\nlimsup
f^n(x^n) \leq f(x)$ by \eqref{eqn:hypoconv}. Thus, the mapping $(x,f)\mapsto f(x)$ on $S\times F$ is
usc and thus measurable. We therefore  have that $((x,y),f)\mapsto
(y-f(x))^2$ is measurable too as a function on $S\times \reals\times F$ and
also a random lsc function. It is trivially locally inf-integrable by its nonnegativity. Theorem \ref{thm:con}(i) therefore applies and a
cluster point $f^\star$ of $\{\hat f^n, n\in\nats\}$ must satisfy a.s.
\[
f^\star \in \nargmin_{f\in F} \Ex\big[(Y^1 - f(X^1)^2\big] = \nargmin_{f\in F} L_P^2(f^0,f)
\]
because $\Ex[Z^1]=0$ and $X^1$ and $Z^1$ are independent; the finite variance of $Z^1$ prevents $\Ex\big[(Y^1 - f(X^1)^2\big]$ from being $\infty$ when $L_P^2(f^0,f)$ is finite. The existence of a cluster point is realized as follows. Let $\phi(f) = \Ex[(Y^1 - f(X^1)^2]$, $f\in F$. If $\{f^n\in F, n\in\nats\}$ satisfies $\setd(f^n,g)\to \infty$ for some $g\in F$, then $\hypo f^n$ set-converges to $\emptyset$ and $f^n(x)\to -\infty$ for all $x\in S$. Thus, $\phi(f^n) = \Ex[(f^0(X^1) - f^n(X^1)^2] + \Ex[(Z^1)^2]\to \infty$ since $f^0$ is real-valued. This implies that $\{f\in F~|~\phi(f) \leq \alpha\}$ is bounded for all $\alpha\in\reals$ and contains $\nOutLim \{f\in F~|~\phi^n(f) \leq \alpha\}$ for any sequence of functions $\{\phi^n:F\to \Reals, n\in\nats\}$ epi-converging to $\phi$ \cite[Thm. 3.1]{BeerRockafellarWets.92} and also $\nlimsup ( \inf_{f\in F} \phi^n(f)) \leq \inf_{f\in F} \phi(f)$ \cite[Thm. 3.8]{RoysetWets.16a}. Under the assumption that $\inf_{f\in F} \phi(f)<\infty$, we therefore have that for some $\bar n$, $\{\epsilon^n\mbox{-}\nargmin_{f\in F} \phi^n(f), n\geq \bar n\}$ is bounded. Applying these facts to the (random) function defined by $\phi^n(f) = n^{-1} \sum_{j=1}^n (Y^j - f(X^j))^2$, which epi-converges to $\phi$ almost surely (cf. Theorem \ref{thm:con}), establishes that $\{\hat f^n, n\in\nats\}$ is bounded almost surely and therefore must have a cluster point. The final
conclusion follows directly from the properties of the $L_P^2$ distance.\eop

\state Proof of Theorem \ref{thm:con_approx}. Following the arguments in the
proof of Theorem \ref{thm:con}, we established that $f\mapsto \phi^n(f) =
n^{-1}\sum_{j=1}^n \psi(X^j,\cdot) + \pi^n$ epi-converges to $f\mapsto
\phi(f) = \Ex[\psi(X^1,\cdot)]$ a.s. as functions on $(F,\setd)$. Next,
suppose that
\[
f^\star \in \nOutLim \big(\epsilon^n\mbox{-}\nargmin_{f\in F^n_\delta} \phi^n(f)\big).
\]
Then there exist a subsequence $\{n_k, k\in \nats\}$ and
\[
f^k \in \epsilon^{n_k}\mbox{-}\nargmin_{f\in F^{n_k}_\delta} \phi^{n_k}(f) \to f^\star.
\]
The continuity of the point-to-set distance and the fact that $\dist(f^k,
F^{n_k}) \leq \delta$ for all $k$ implies that $\dist(f^\star,\nLim F^n) \leq
\delta$, i.e., $f^\star \in F_\delta^\infty$. Thus, it only remains to show
that $\phi(f^\star) \leq \ninf_{f\in \nLim F^n} \phi(f)$. Let $g^\star \in \nargmin_{f\in \nLim F^n} \phi(f)$. Then, because $\phi^n$ epi-converges to $\phi$, there exists $g^n\in F\to
g^\star$ such that
\[
\nlimsup \phi^n(g^n) \leq \phi(g^\star).
\]
Since $g^\star \in \nLim F^n$, there is $\bar n\in\nats$ such that
$\dist(g^n,F^n) \leq \delta$ for all $n\geq \bar n$. Consequently, leveraging
the epi-convergence property and the above facts,
\begin{align*}
  \phi(f^\star) &\leq \nliminf \phi^{n_k}(f^{k}) \leq \nliminf \Big(\ninf_{f\in F_\delta^{n_k}} \phi^{n_k}(f) + \epsilon^{n_k}\Big)\\
                         &\leq \nlimsup \phi^{n_k}(g^{n_k}) \leq \phi(g^\star) = \ninf_{f\in \nLim F^n} \phi(f).
\end{align*}
The first conclusion is established. The second conclusion is immediate after
realizing that $F_\delta^\infty = F$ when $\nLim F^n = F$.\eop

\state Proof of Corollary \ref {cor:conML_approx}. The arguments of Corollary
\ref{cor:conML} in conjunction with Theorem \ref{thm:con_approx} yield
$f^\star \in F_\delta^\infty$ and $K(f^0; f^\star) \leq \ninf_{g\in \nLim
F^n} K(f^0; g)$. Since $\nLim F^n \subset F$ consists only of densities and
$f^0 \in \nLim F^n$, the right-hand side in this inequality is zero and the
conclusion follows.\eop

\state Proof of Theorem \ref{thm:plug}. The assertions about $\hat m^n$ and
$\hat h^n$ are essentially in \cite[Prop. 2.1]{Royset.16}, with an
extension to improper functions following straightforwardly. The conclusion about $\hat l^n$ holds by \cite[Prop. 7.7]{VaAn}.\eop

\begin{lemma}{\rm (hypo-convergence under composition).}\label{lem:composition}
For a continuous nondecreasing function $h_0:\reals\to \Reals$, let $h:\Reals \to \Reals$ have $h(y) = h_0(y)$ if $y\in\reals$, $h(-\infty)=\inf_{\bar y\in\reals} h_0(\bar y)$, and $h(\infty)=\sup_{\bar y\in\reals} h_0(\bar y)$. If $g^n:S\to \Reals$ hypo-converges to $g:S\to \Reals$, then $h\comp g^n$ hypo-converges to $h\comp g$.
\end{lemma}
\state Proof. Suppose that $x^n\in S\to x$, which implies that $\nlimsup g^n(x^n) \leq g(x)$. Fix $n$ and let $\epsilon>0$. Suppose that $\xi^n = \sup_{m\geq n} h(g^m(x^m))\in\reals$. Then, there exists $\bar m\geq n$ such that $\xi^n \leq h(g^{\bar m}(x^{\bar m})) + \epsilon \leq h(\sup_{m\geq n} g^m(x^m)) + \epsilon$, the last inequality holds because $h$ is nondecreasing. Since $\epsilon$ is arbitrary, $\xi^n\leq h(\sup_{m\geq n} g^m(x^m))$. A similar argument leads to the same inequality if $\xi^n = \infty$ and, trivially, also when $\xi^n = -\infty$.
Since the inequality holds for all $n$, it follows by the continuity of $h$ that
\begin{align*}
\nlimsup h(g^n(x^n)) & = \lim_{n\to\infty} \big(\nsup_{m\geq n} h(g^m(x^m))\big) \leq \lim_{n\to\infty} h\big(\nsup_{m\geq n} g^m(x^m)\big)\\
& = h\big( \lim_{n\to\infty} ( \nsup_{m\geq n} g^m(x^m))\big) = h(\nlimsup g^n(x^n)) \leq h(g(x)).
\end{align*}
For any $x\in S$, there exists $x^n\in S\to x$ with $g^n(x^n) \to g(x)$. Since $h$ is continuous, this implies  $h(g^n(x^n)) \to h(g(x))$ and the conclusion follows.\eop

\state Proof of Proposition \ref{prop:closedness_convex}. The first claim follows by \cite[Prop. 4.15]{VaAn}. Since $-f^n,
-f$ are proper, lsc, and convex, it follows by \cite[Thm. 12.35]{VaAn}
that the graphs of the subdifferentials $\partial f^n$ set-converge to the
graph of $\partial f$. Thus, for every $(x,v)$ in the graph of $\partial f$,
there exists $x^n\to x$ and $v^n\to v$, with $v^n\in
\partial f^n(x^n)$. Since $\|v^n\|_2\leq \kappa$ for all $n$, we also have
that $\|v\|_2\leq \kappa$.

For part (ii), Lemma \ref{lem:composition}, with $h$ defined by $h(y) = \log y$ if $y\in (0,\infty)$, $h(y) = -\infty$ if $y = [-\infty,0]$, and $h(y) = \infty$ if $y = \infty$, yields that $h\comp f^n$ hypo-converges to $h\comp f$. Since $h\comp f^n$ is concave, it follows by \cite[Prop. 4.15]{VaAn} that $h\comp f$ is concave too.

For part (iii), let $\lambda \in (0,1)$ and $x,y\in \nt S$. Set $z = \lambda
x + (1-\lambda)y$. Hypo-convergence implies that there exists $z^n\in \nt
S\to z$ such that $f^n(z^n) \to f(z)$. Construct $x^n = x + z^n - z$ and $y^n
= y + z^n - z$. Clearly, $x^n\to x$ and $y^n\to y$. Then, $\lambda x^n +
(1-\lambda) y^n = z^n$. Let $\epsilon>0$ and suppose that $f(z)<\infty$,
$f(x)>-\infty$, and $f(y) > -\infty$. There exists $\bar n$ such that for all
$n\geq \bar n$, $x^n, y^n\in S$ and
\begin{equation*}
f(z) \leq f^n(z^n) + \frac{\epsilon}{3}, ~~ f^n(x^n) \leq f(x) + \frac{\epsilon}{3\lambda}, ~~ f^n(y^n) \leq f(y) + \frac{\epsilon}{3(1-\lambda)}.
\end{equation*}
Collecting these results and use the convexity of $f^n$, we obtain that for
$n\geq \bar n$
\begin{align*}
  f(z) & \leq f^n(z^n) + \frac{\epsilon}{3}\leq \lambda f^n(x^n) + (1-\lambda) f^n(y^n) + \frac{\epsilon}{3}\\
  & \leq  \lambda f(x) +  (1-\lambda)f(y) + \epsilon.
\end{align*}
Since $\epsilon>0$ is arbitrary, $f(z) \leq \lambda f(x) + (1-\lambda)f(y)$.
A similar argument leads to the same conclusion when $f(z)=\infty$,
$f(x)=-\infty$, and/or $f(y) = -\infty$.

It only remains to examine the case when $x$ and/or $y$ are at the boundary
of $S$. Suppose that $\lambda\in (0,1)$, $x\in \nt S$, and $y \in
S\setminus\nt S$. Then, there exists $y^n\in\nt S\to y$ with $f(\lambda x +
(1-\lambda) y^n) \leq \lambda f(x) + (1-\lambda) f(y^n)$ because $S$ must be
convex. Since $\lambda x + (1-\lambda) y^n, \lambda x + (1-\lambda) y\in \nt
S$ and $f$ is continuous on $\nt S$, the left-hand side tends to $f(\lambda x
+ (1-\lambda) y)$. The upper limit of the right-hand side is   $\lambda f(x)
+ (1-\lambda) f(y)$ by the usc of $f$.  A similar argument holds in the other
cases. Thus, $f$ is convex.\eop

\state Proof of Proposition \ref{prop:closedness_trans}. By \cite[Thm. 7.6]{VaAn}, either $\hypo g^n$ set-converges to $\emptyset$ or there exist $g\in \uscfcns(S)$ and a subsequence $\{n_k, k\in\nats\}$ such that $g^{n_k}\to g$. In the second case, $h\comp g^{n_k}$ hypo-converges to $h\comp g$ by Lemma \ref{lem:composition}. Since a hypo-limit is unique, $f = h\comp g$. In the first case, for all $x\in S$, $g^n(x)\to-\infty$ so that $h\comp g^n(x) \to f(x) = h(-\infty) = h\comp g(x)$, when $g(x) = -\infty$ for all $x\in S$. The conclusions then follow by Proposition \ref{prop:closedness_convex}.\eop

\state Proof of Proposition \ref{prop:closedness_mono}. For part
(i), let $x\leq y$, with $y\in \nt S$, and $\epsilon>0$. The usc property
implies that there exists $\delta>0$ such that $f(y) \geq f(z)-\epsilon$ for
all $z\in S$ with $\|z-y\|\leq \delta$. Since $y\in \nt S$, $z$ can be
takes such that $z_i > y_i$ for $i=1, \dots, d$ and $z\in \nt S$. By
hypo-convergence, there exists $x^n\in S\to x$ such that $f(x) \leq \nliminf
f^n(x^n)$ and also $\nlimsup f^n(z) \leq f(z)$. Thus, $x^n \leq z$ for
sufficiently large $n$. By the nondecreasing property,
\begin{equation*}
  f(x) \leq \nliminf f^n(x^n) \leq \nliminf f^n(z) \leq \nlimsup f^n(z) \leq f(z)  \leq f(y) + \epsilon.
\end{equation*}
Since $\epsilon>0$ is arbitrary the first conclusion follows.

Under the additional structure of $S$, the argument can be modified as
follows. Now with $y\in S$, let $\delta>0$ and $x^n$ be as earlier. Construct
$z\in \reals^d$ by setting $z_i = \min\{\beta_i, y_i + \delta\}$. Let $\bar
n$ be such that $x_i^n \leq x_i + \delta$ for all $i=1, \dots, d$ and $n\geq
\bar n$. Then, for $n\geq\bar n$, $x_i^n$ $\leq$ $\min\{\beta_i, x_i +
\delta\}$ $\leq$ $\min\{\beta_i, y_i + \delta\}= z_i$. Thus, again we have
that $x^n \leq z$ for sufficiently large $n$ and the preceding arguments lead
to the conclusion.

For (ii) let $x\leq y$, with $x\in \nt S$, and $\epsilon>0$. The usc property
implies that there exists $\delta>0$ such that $f(x) \geq f(z)-\epsilon$ for
all $z\in S$ with $\|z-x\|\leq \delta$. Since $x\in \nt S$, $z$ can be
takes such that $z_i < x_i$ for $i=1, \dots, d$ and $z\in \nt S$. In view of
the hypo-convergence, there exists $y^n\in S\to y$ such that $f(y) \leq
\nliminf f^n(y^n)$ and also $\nlimsup f^n(z) \leq f(z)$. Thus, $z\leq y^n$
for sufficiently large $n$. Using the nonincreasing property, we then obtain
that
\begin{equation*}
  f(y) \leq \nliminf f^n(y^n) \leq \nliminf f^n(z) \leq \nlimsup f^n(z) \leq f(z)  \leq f(x) + \epsilon.
\end{equation*}
Since $\epsilon>0$ is arbitrary the first conclusion follows.

Under the additional structure of $S$, the argument can be modified as
follows. Now with $x\in S$, let $\delta>0$ and $y^n$ be as earlier. Construct
$z\in \reals^d$ be setting $z_i = \max\{\alpha_i, x_i - \delta\}$. Let $\bar
n$ be such that $y_i^n \geq y_i - \delta$ for all $i=1, \dots, d$ and $n\geq
\bar n$. Then, for $n\geq\bar n$, $y_i^n$ $\geq$ $\max\{\alpha_i, y_i -
\delta\}$ $\geq$ $\max\{\alpha_i, x_i - \delta\} = z_i$. Again we have $z
\leq y^n$ for sufficiently large $n$ and the preceding arguments lead to the
conclusion.\eop

\state Proof of Proposition \ref{prop:closedness_lip}. If $\kappa =
0$, then $f^n$ are constant functions on $S$ and $f$ also, and the conclusion
holds. Suppose that $\kappa>0$. Let $x,y\in S$, with $f(x)$ and $f(y)$
finite,  and $\epsilon>0$. Hypo-convergence implies that there exists $x^n\in
S\to x$ such that $f^n(x^n) \to f(x)$ and $\nlimsup f^n(y) \leq f(y)$. Hence,
there exists $\bar n$ such that for all $n\geq \bar n$, $\|x^n-x\|\leq
\epsilon/(3\kappa)$, $|f^n(x^n) - f(x)| \leq \epsilon/3$, $f^n(y) \leq f(y) +
\epsilon/3$. For such $n$, $f(x) - f(y)$
\begin{align*}
   & = f(x) - f^n(x^n) + f^n(x^n) - f^n(x) + f^n(x) - f^n(y) + f^n(y) - f(y)\\
   & \leq \frac{\epsilon}{3} + \kappa\|x^n - x\| + \kappa\|x-y\| + f(y) + \frac{\epsilon}{3} - f(y) \leq \kappa\|x-y\| + \epsilon.
\end{align*}
Repeating this argument with the roles of $x$ and $y$ interchanged, we obtain
that $|f(x) - f(y)| \leq \kappa\|x-y\| + \epsilon$. Since $\epsilon>0$ is
arbitrary, $f$ is Lipschitz continuous with modulus $\kappa$ when finite. If
$f$ is not finite on $S$, then it cannot be Lipschitz continuous.\eop

\state Proof of Proposition \ref{prop:closedness_bd}. Let $x\in S$
and observe that $g(x) \leq \nlimsup f^n(x)$ $\leq$ $f(x)$ by \eqref{eqn:hypoconv}, which established
the lower bound. Since $h$ is usc, we also have that for some $x^n\in S\to
x$, $h(x) \geq \nlimsup h(x^n)$ $\geq$ $\nliminf f^n(x^n)$ $\geq f(x)$, which
confirms the upper bound.\eop

\state Proof of Proposition \ref{prop:mtp2}. Since the collection of
functions is equi-usc, hypo-convergence implies pointwise convergence (Proposition \ref{prop:equi-usc_ptwise}) and the conclusion follows immediately.\eop

\state Proof of Proposition \ref{prop:lscsupnorm}. Let $\epsilon>0$ and $f^n\in F\to f$. First, suppose that $\nsup_{x\in S} g(f(x))\in \reals$. Then, there exists $\bar x\in S$ such that $g(f(\bar x)) \geq \nsup_{x\in S} g(f(x)) - \epsilon$. By \eqref{eqn:hypoconv}, there is $x^n\in S\to \bar x$ such that $f^n(x^n)\to f(\bar x)$. Since $g$ is lsc,
\[
\nliminf \big( \nsup_{x\in S} g(f^n(x))\big) \geq \nliminf g(f^n(x^n)) \geq g(f(\bar x)) \geq \nsup_{x\in S} g(f(x)) - \epsilon.
\]
Second, suppose that $\nsup_{x\in S} g(f(x)) = \infty$. Then, there exists $\bar x\in S$ such that $g(f(\bar x)) \geq 1/\epsilon$. Again, there is $x^n\in S\to \bar x$ such that $f^n(x^n)\to f(\bar x)$ and
\[
\nliminf \big( \nsup_{x\in S} g(f^n(x))\big) \geq \nliminf g(f^n(x^n)) \geq g(f(\bar x)) \geq 1/\epsilon.
\]
Since $\epsilon>0$ is arbitrary, we have established that $\nliminf \big( \nsup_{x\in S} g(f^n(x))\big) \geq \nsup_{x\in S} g(f(x))$; it trivially holds when $\nsup_{x\in S} g(f(x))$ $= -\infty$.\eop

\state Proof of Proposition \ref{prop:integral2}. Since the collection of
functions is equi-usc at Lebesgue-a.e. $x\in S$, hypo-convergence implies pointwise convergence at Lebesgue-a.e. $x\in S$ by Proposition \ref{prop:equi-usc_ptwise}. The
conclusions follow directly from an application of the Dominated Convergence
Theorem.\eop

\state Proof of Proposition \ref{prop:moment}. Since $C\subset C^n$, $F\subset F^n$ and it suffices to confirm
that $\nOutLim F^n \subset F$. Take $f\in \nOutLim F^n$. There exists $f^k\in
F^{n_k}\to f$. Since $\int x f^k(x) dx \in C^{n_k}$, that integral
converges to $\int x f(x) dx$ by Proposition \ref{prop:integral2}, and the set-convergence $C^n$ to $C$ allow us to conclude that $\int xf(x)
dx \in C$.\eop

\state Proof of Theorem \ref{thm:algo}. Let $\phi,\phi^\nu:F^0\to
(-\infty,\infty]$ be given by $\phi(f)  = n^{-1}\sum_{j=1}^n \psi(x^j,f) +
\pi(f)$ if $f\in F$ and $\phi(f) = \infty$ otherwise; and $\phi^\nu(f)  =
n^{-1}\sum_{j=1}^n \psi(x^j,f) + \pi^\nu(f)$ if $f\in F^\nu$ and $\phi^\nu(f)
= \infty$ otherwise. We start by showing that $\phi^\nu$ epi-converges to
$\phi$. Let $f^\nu \in F^0\to f$. If $f\in F$, then
\[
\nliminf \phi^\nu(f^\nu) \geq \frac{1}{n}\sum_{j=1}^n \psi(x^j,f) + \pi(f) = \phi(f).
\]
If $f\not\in F$, then because $F$ is closed we must have that $f^\nu\not\in
F^\nu$ for sufficiently large $\nu$. Thus, $\nliminf \phi^\nu(f^\nu) =
\phi(f)=\infty$. Next, let $f\in F$. There exists $f^\nu\in F^\nu\to f$ because
$F^\nu$ set-converges to $F$. Then,
\[
\nlimsup \phi^\nu(f^\nu) = \nlimsup \Big(\frac{1}{n}\sum_{j=1}^n \psi(x^j,f^\nu) + \pi^\nu(f^\nu) \Big) = \phi(f).
\]
This is sufficient for $\phi^\nu$ epi-converging to $\phi$. Reasoning along
the lines of those in the proof of Theorem \ref{thm:con} yields (i).

For (ii), we recognize that the additional condition on $F^0$ ensures that it
is compact. Thus, $\{f^\nu, \nu\in\nats\}$ in the statement of the theorem
must have a cluster point. Every such cluster point must be in
$\epsilon^\infty\mbox{-}\nargmin_{f\in F} \phi(f)$. Let $\delta =
\epsilon-\epsilon^\infty$, which is positive. Since $F^0$ is compact,
 $\pi^\nu$ converges uniformly
to $\pi$. Hence, there exists $\bar \nu\in\nats$ such that $\pi(f^\nu) \leq
\pi^\nu(f^\nu) + \delta/3$, $\epsilon^\nu \leq \epsilon^\infty + \delta/3$,
and, in view of epi-convergence, $\ninf_{f\in F^\nu} \phi^\nu(f) \leq
\ninf_{f\in F} \phi(f) + \delta/3$ for all $\nu\geq \bar \nu$. Since
$F^\nu\subset F$, we then have
\begin{align*}
 & \phi(f^{\bar \nu})  = \frac{1}{n}\sum_{j=1}^n \psi(x^j,f^{\bar \nu}) + \pi(f^{\bar \nu}) \leq \frac{1}{n}\sum_{j=1}^n \psi(x^j,f^{\bar \nu}) + \pi^{\bar \nu}(f^{\bar \nu}) + \delta/3\\
    \leq &\ninf_{f\in F^{\bar\nu}} \phi^{\bar \nu}(f) + \epsilon^{\bar \nu} + \delta/3 \leq \ninf_{f\in F} \phi(f) + \epsilon^\infty + \delta = \ninf_{f\in F} \phi(f) + \epsilon,
\end{align*}
which establishes the claim.\eop

\end{document}